\documentclass[a4paper,12pt]{article}

\usepackage{amsmath,amsfonts,color,enumitem,geometry,latexsym,mathrsfs,upgreek}
\usepackage[colorlinks,linkcolor=blue,citecolor=red,urlcolor=blue]{hyperref}
\usepackage[amsmath,thmmarks]{ntheorem}
\usepackage{graphicx}
\usepackage{cleveref}
\usepackage{nameref}

\newcommand{\Title}{On equilibrium states for certain partially hyperbolic endomorphisms with one-dimensional center}
\newcommand{\Titleabbr}{On equilibrium states for certain partially hyperbolic endomorphisms with one-dimensional center}
\date{}
\newcommand{\Keyword}{partially hyperbolic endomorphism; equilibrium state; unstable pressure; stable pressure; inverse limit space}
\newcommand{\Fund}{This work is supported by  NSFC (Noes: 12171400, 12426107, 12571206).}
\newcommand{\Authora}{Yifan Zhang${}^\text{ *}$}
\newcommand{\AuthorabbrA}{ZHANG}
\newcommand{\Authorb}{Yujun Zhu}
\newcommand{\AuthorabbrB}{ZHU}
\newcommand{\Emaila}{exampleA@email.com}
\newcommand{\Emailb}{exampleB@email.com}

\usepackage{fancyhdr}
\title{\Title
	\footnotetext{
		\setlength{\parindent}{2em}
		\emph{2020 Mathematics Subject Classification: {37D30, 37D35.}}\\
		\indent \emph{Keywords and phrases:} \Keyword.\\
		\indent \Fund\\
		\indent{} $^*$ The corresponding author.}}
\author{\Authora \text{ and}  \Authorb }
\date{}
{
	\theorembodyfont{\rm}
	\newtheorem{mdef}{Definition}[section]
	\newtheorem{rmk}{Remark}
	\newtheorem{exa}{Example}
}
\newtheorem{thm}{Theorem}[section]

\newtheorem*{mainthm}{Main Theorem}
\newcommand{\thmref}[1]{\hyperref[#1]{Main Theorem}}
\newtheorem{prop}[thm]{Proposition}
\newtheorem{lem}[thm]{Lemma}
\newtheorem{cor}[thm]{Corollary}
\newtheorem{claim}[thm]{Claim}

\crefname{section}{§}{§§}
\Crefname{section}{§}{§§}

\theoremstyle{nonumberchange}
\theorembodyfont{\rm}
\theoremsymbol{\mbox{$\Box$}}
\newtheorem{pf}{Proof}
\newtheorem{pot3.1}{Proof of Proposition {\ref{3.1}}}
\newtheorem{pot4.3}{Proof of Proposition {\ref{4.3}}}

\begin{document}
	\vskip 5.1cm
	
	\maketitle		
	
	\begin{center}
		\small
		\textsl{School of Mathematical Sciences, \\Xiamen University, Xiamen, 361005, PR~China}
	\end{center}
	\begin{center}
		\begin{minipage}{130mm}
			
			{\bf Abstract}:
			In this paper, the equilibrium states for a non-degenerate $ C^2 $ partially hyperbolic endomorphism $f$ on a closed Riemannian manifold $M$ with one-dimensional center bundle are investigated. Applying the criterion of Climenhaga-Thompson (\cite{CT21}) and the method of Mongez-Pacifico (\cite{Mongez}), we use the techniques of inverse limit to obtain the uniqueness and robustness of equilibrium states for $f$ and any H\"{o}lder continuous potential satisfying certain conditions about the unstable pressure and stable pressure.
		\end{minipage}
		
	\end{center}
	
	\section{Introduction}
	In modern dynamical systems and ergodic theory, both entropy and pressure are the fundamental quantities that measure the complexity of a system. These concepts are closely related to physics and thermodynamic mechanisms; see the monographs by Ruelle \cite{Ruelle} and Bowen \cite{Bowen book}. Roughly speaking, the topological pressure can be interpreted as the maximal free energy that a system can attain. Then, the Gibbs distribution, which can maximize free energy, corresponds to the equilibrium state in ergodic theory. If one neglects the potential energy, topological pressure reduces to topological entropy.
	
	Invariant measures play central roles in ergodic theory. Let $f$ be a continuous map on a compact metric space $X$. Given an $f$-invariant measure $\mu$, one can define the measure-theoretic entropy $h_{\mu}(f)$. The variational principle establishes a relationship between the measure-theoretic entropy and the topological entropy  $h(f)$ as
	$$
	h(f)=\sup\{h_{\mu}(f) \;|\; \mu\in \mathcal{M}_f\},
	$$
	where $ \mathcal{M}_f $ is the set of $ f $-invariant Borel probability measures.	The element in $ \mathcal{M}_f $ that makes the equality hold is called the \emph{measure of maximal entropy}.
	Given a continuous function $\phi$ on $X$, referred to as a potential, the variational principle for pressure then yields
	$$
	P(f,\phi)=\sup\{h_{\mu}(f)+\int \phi d\mu \;|\; \mu\in \mathcal{M}_f\}.
	$$
	The element in $ \mathcal{M}_f $ that reaches the supremum is called the \emph{equilibrium state} for $\phi$. It follows that when the potential function vanishes, the topological pressure reduces to the topological entropy, and the corresponding equilibrium state coincides with a measure of maximal entropy.
	
	The investigations into the existence and uniqueness of the equilibrium state and the measure of maximal entropy are of interest to the researchers. Among the techniques in this topic,  there are two classical ones established in the  1970s. The first one is using the transfer operator theory, see in \cite{Bowen book} and \cite{Ruelle}. For a topologically transitive locally maximal hyperbolic set $\Lambda$ of a diffeomorphism $f$, or especially, a basic set of an Axiom A diffeomorphism, one can get a Markov partition of $\Lambda$ and then obtain a symbolic extension of $f_\Lambda$. Hence, applying the transfer operator theory for a finite-type symbolic system, one obtains the existence and uniqueness of the equilibrium state for certain potentials. The second one is using the expansivity and the specification property, see Bowen \cite{Bowen74}. Precisely, for a homeomorphism $f$ on a compact metric space which is expansive and has the specification property, Bowen showed that for a potential function satisfying certain conditions (called the Bowen property), there exists a unique equilibrium state. In fact, if $\Lambda$ is a topologically mixing compact locally maximal hyperbolic set
	for a diffeomorphism $f$, then $f_\Lambda$ is expansive and has the specification property. Therefore, for such type of uniformly hyperbolic systems, there are two ways to investigate the existence and uniqueness of the equilibrium state and the measure of maximal entropy.
	
	In recent years, the study of dynamics of the differentiable systems beyond uniform hyperbolicity, such as non-uniformly hyperbolic systems and partially hyperbolic systems,  has attracted much attention and made a series of important advances. For partially hyperbolic diffeomorphisms, a center direction is allowed in addition to the hyperbolic directions. The presence of the center direction permits a very rich type of structure. Investigating the thermodynamic mechanisms, especially the existence and uniqueness of the equilibrium state, for partially hyperbolic systems require more exquisite techniques.
	
	Recently, Climenhaga and Thompson \cite{CT21} ingeniously developed the criterion of Bowen \cite{Bowen74} in settings beyond uniform hyperbolicity
	using weakened versions of specification and expansivity. They first considered the case of symbolic systems.
	Rather than assuming specification on the entire symbolic system, they identified a collection of ``good" words on which the specification property holds, while showing that the set of ``bad" words has topological entropy strictly smaller than that of the full system. Using this decomposition, they established the uniqueness of the measure of maximal entropy and the equilibrium state in \cite{CT12} and \cite{CT13}, respectively. Further, Climenhaga and Thompson weakened expansivity by introducing a weaker, non-uniform notion applicable to general compact metric spaces; they also formulate correspondingly weakened, non-uniform versions of specification, expansivity, and the Bowen property in \cite{CT14} and \cite{CT16}. For applications of the method mentioned above, one can see \cite{CFT18} and \cite{CFT19}. The general strategy is to first identify obstructions to expansivity and specification, and then to control their entropy or pressure. Pacifico, F. Yang, and J. Yang \cite{Yang} further refined the conditions in the Climenhaga-Thompson approach, allowing the framework to be applied in settings with fixed parameters.
	
	From \cite{CT21}, one can observe that for a partially hyperbolic diffeomorphism with one-dimensional center, if every leaf of the stable and unstable foliations is dense, and the supremum of the measure-theoretic entropy over ergodic measures with non-positive or non-negative central Lyapunov exponent is different, then there exists a unique measure of maximal entropy. A similar argument applies to the case of pressure. Mongez and Pacifico \cite{Mongez} substituted the condition on the central Lyapunov exponent used in \cite{CT21} with the requirement that there exists a strict gap between the unstable entropy and the stable entropy of $ f $. The definitions of unstable entropy and stable entropy referenced here are drawn from \cite{Yan21, HHW}.  Liu and Liao  \cite{LL} showed the robustness and uniqueness of equilibrium states for a class of partially hyperbolic with a central direction that admits a dominated splitting into one-dimensional subbundles and H\"{o}lder continuous potentials with not very large oscillation.
	
	Partially hyperbolic endomorphisms can be viewed as a generalization of partially hyperbolic diffeomorphisms. However, several distinctions arise when attempting to establish analogous results in this broader setting. Due to the non-invertibility of the map, each point may possess multiple preimages, which prevents the direct definition of unstable manifolds in the classical sense. To address this issue, we adopt the framework of inverse limit spaces for partially hyperbolic endomorphisms, building on the foundational work presented in \cite{Qian}. Through the inverse limit space, the unstable manifolds of partially hyperbolic endomorphisms can be rigorously defined, and corresponding notions of unstable entropy and stable entropy have also been introduced in \cite{Wang, Zhu}.   Álvarez and Cantarino \cite{Cantarino} established the existence of equilibrium states for partially hyperbolic endomorphisms with one-dimensional center bundles by showing that the shift map on the inverse limit space is $h$-expansive, as well as the uniqueness of such states for certain endomorphisms defined on the torus. Recently, Arbieto and  Cabezas  \cite{AC} established the existence of equilibrium states for a class of partially hyperbolic endomorphisms when the central direction is simple or which is decomposed into one-dimensional subbundles with dominated splitting, and gave the finiteness of measures of maximal entropy for a certain partially hyperbolic endomorphism with one-dimensional central direction.
	
	Our main purpose in this paper is to generalize the results in \cite{Mongez} to the case of partially hyperbolic endomorphisms. Here we assume that the map $ f $ is non-degenerate and thus $ f $ has uniform separation of preimages (see Lemma \ref{separation}), then we can deal with the preimages on a small domain. Denote  $U_a(\phi):=\{\psi\in C(X)\;|\; \lVert \phi-\psi \rVert <a\}$, where $a>0 $ and $C(X)$ represents all continuous functions on $X$ and $ \lVert \phi \rVert=\sup_{x\in X}\lvert \phi(x)\rvert $.
	
	\begin{mainthm}\label{Main}
		Let $ f:M\to M $ be a non-degenerate $ C^2 $  partially hyperbolic endomorphism on a compact Riemannian manifold $ M $ with one-dimensional center bundle and $ \phi:M\to \mathbb{R} $ be a H\"{o}lder continuous potential. Assume that  $ P^u(f,\phi)>P^s(f,\phi)$ and  the unstable foliation $\tilde{ \mathcal{W}}^u $ of $f$ is minimal, then for any $ 0<a<\frac{P^u(f,\phi)-P^s(f,\phi)}{2} $
		there exists a $ C^1 $ neighborhood $ \mathcal{U}_f $ of $ f $ such that for any H\"{o}lder continuous potential $\psi\in U_{a}(\phi)$ and each non-degenerate $ C^2 $ partially hyperbolic endomorphism $ g\in \mathcal{U}_f $ there   has a unique equilibrium state for $(g,\psi)$.
	\end{mainthm}
	
	\begin{rmk}
		The results of the Main Theorem still hold for $ C^{1+\alpha} $ mappings. To cite some conclusions in \cite{Qian} conveniently, we therefore assume $ C^2 $ regularity here.
	\end{rmk}
	
	The proof of the main theorem in this paper is inspired by arguments established in \cite{CT21} and \cite{Mongez}. Compared with the case of diffeomorphisms, the main difficulty here comes from non-invertibility. In order to overcome this deficiency, we use the technique of inverse limit. For certain non-degenerate partially hyperbolic endomorphisms with one-dimensional center and any  H\"{o}lder continuous potential function satisfy specific conditions, we aim to construct a decomposition of orbit segments in the inverse limit spaces. An important condition in the main theorem is the gap between the unstable pressure and the stable pressure. For partially hyperbolic diffeomorphisms, the unstable entropy and unstable pressure for a continuous function were introduced in \cite{HHW, HWZ}, respectively. As a generalization, the concepts of unstable entropy and unstable pressure were given by Wang, Wu and Zhu in \cite{Wang} for partially hyperbolic endomorphisms. In \cite{Zhu}, Wu and Zhu introduced the stable entropy for non-degenerate $C^1$  partially hyperbolic endomorphisms and their stable pressure is given in \cite{Li}. Different from the condition given by Mongez and Pacifico in \cite{Mongez}, we change the condition ``$h^u(f)-h^s(f)>\sup\phi-\inf\phi\geq 0$'' to ``$P^u(f,\phi)>P^s(f,\phi)$'', where
	$P^u(f,\phi)=\sup\{h^u_{\mu}(f)+\int \phi d\mu\;|\; \mu\in \mathcal{M}_f^e\}$ and $P^s(f,\phi)=\sup\{h^s_{\mu}(f)+\int \phi d\mu\;|\; \mu\in \mathcal{M}_f^e\}$.
	Since for each ergodic $f$-invariant measure $\mu$, there is
	\begin{align*}
		h^u_\mu(f)+\min\phi \leq h^u_\mu(f)+\int \phi d\mu,
	\end{align*}
	we have $P^u(f,\phi)\geq h^u(f)+\min\phi$.
	Meanwhile, $h^s(f)+\max\phi \geq P^s(f,\phi)$. Therefore, the condition ``$P^u(f,\phi)>P^s(f,\phi)$'' is weaker than ``$h^u(f)-h^s(f)>\sup\phi-\inf\phi\geq 0$'' since the space we are considering is compact so $\sup\phi = \max\phi$ and $\inf\phi=\min\phi$.
	Under the condition $ P^u(f,\phi)>P^s(f,\phi) $ and assuming the existence of unstable foliation with some sense of minimality, we show that the specification property holds for good orbit segments that exhibit a contracting property in the center-stable direction. Moreover, the potential $ \phi $ satisfies the Bowen property on these good orbit segments. Additionally, the topological pressure on the collection of bad orbit segments and the pressure associated with the obstruction to expansivity are both strictly less than the topological pressure of the entire system. As a consequence, there exists a unique equilibrium state on the inverse limit space, and its projection yields the unique equilibrium state for the original partially hyperbolic endomorphism. Furthermore, we also demonstrate that the conditions of the theorem are inherently robust, which further implies the robustness of the equilibrium state.
	At the same time, H\"{o}lder continuous potential functions can also vary within a certain range,  thus allowing the uniqueness of the equilibrium state of the system under small perturbation can be obtained.
	
	This paper is organized as follows. In Section 2, we give some preliminaries, including the concepts of the decomposition of the orbit segments, the pressure of the obstructions to expansivity, partially hyperbolic endomorphisms, inverse limit space, (un)stable en        tropy and (un)stable pressure, and other necessary definitions. In Section 3, we demonstrate the robustness of the conditions in Main Theorem. In Section 4, we prove  Main Theorem using the improved Climenhaga-Thompson criterion. In Section 5, we give an example satisfying the conditions in our main theorem.

	\section{Preliminaries and fundamental properties}
	\subsection{Conditions regarding the uniqueness of equilibrium state}
	In this section, we give some notions which will be used to prove the uniqueness of the equilibrium state from \cite{CT12, CT14, CT16, Yang}. Let $ f:X\to X $ be a continuous map on a compact metric space $ X$.
	\subsubsection{\texorpdfstring{Topological pressure}{Topological pressure}}
	Given $ n\in\mathbb{N},\delta>0 $ and $ x\in X $, the Bowen ball of order $n  $   with radius $ \delta $ at $ x $ is defined by
	$$
	B_n(x,\delta):=\{y\in X \;|\; d_n(x,y)<\delta\},
	$$
	where $ d_n(x,y):=\max\{d(f^ix,f^iy)\;|\;0\leq i<n\} $. A subset $ E$ of $ X $ is said to be $ (n,\delta) $-separated with respect to $ f $ if $ x,y\in E, x\neq y $, implies $ d_n(x,y)>\delta $.
	
	We identify $ X\times\mathbb{N} $ with the space of finite orbit segments by associating each pair $ (x,n)\in X\times \mathbb{N} $ with $ \{x,fx,\cdots,f^{n-1}x\} $. Given $ \mathcal{C}\subset X\times \mathbb{N} $ and $ n\in \mathbb{N} $ we denote $ \mathcal{C}_n:=\{x\in X \;|\; (x,n)\in\mathcal{C}\} $.
	
	Let $ \phi:X\to \mathbb{R} $ be continuous, $ n\geq 1 $ and $ \varepsilon>0 $. We denote $ \sup_{y\in B_n(x,\varepsilon)}\sum_{i=0}^{n-1}\phi(f^iy) $ by $ \Phi_{\varepsilon}(x,n) $. In particular, $ \Phi_0(x,n)=\sum_{i=0}^{n-1}\phi(f^ix) $. For  $\varepsilon>0$, $\delta>0$, $\mathcal{C}\subset X\times\mathbb{N}$ and $ n\in \mathbb{N} $, put
	\begin{equation*}
		\Lambda(\mathcal{C},f,\phi,\delta,\varepsilon,n)=\sup\left\{\sum_{x\in E}e^{\Phi_\varepsilon(x,n)}\;|\;E\subset \mathcal{C}_n \text{ is }(n,\delta)\text{-separated}\right\},
	\end{equation*}
	which is said to be the \textit{partition function}. Without causing any ambiguity, we can omit $ f $ in $ \Lambda(\mathcal{C},f,\phi,\delta,\varepsilon,n) $. We denote $ \Lambda(\mathcal{C},\phi,\delta,n):=\Lambda(\mathcal{C},\phi,\delta,0,n) $ whenever $ \varepsilon=0 $.
	
	The pressure of $ \phi  $ on $ \mathcal{C} $ at scale $ \delta,\varepsilon $ is defined by
	\begin{equation*}
		P(\mathcal{C},\phi,\delta,\varepsilon)=\limsup_{n\to \infty}\frac{1}{n}\log \Lambda(\mathcal{C},\phi,\delta,\varepsilon,n).
	\end{equation*}
	We write $ P(\mathcal{C},\phi,\delta) $ instead of $ P(\mathcal{C},\phi,\delta,0) $ whenever $ \varepsilon=0 $. Put
	\begin{equation*}
		P(\mathcal{C},\phi)=\limsup_{\delta\to 0}P(\mathcal{C},\phi,\delta).
	\end{equation*}
	The limit exists since $ P(\mathcal{C},\phi,\delta) $ is monotonic when considered as a function of $ \delta $. Taking $ \mathcal{C}=X\times\mathbb{N} $, we have $ P(f,\phi):=P(\mathcal{C},\phi) $ as the usual notion of topological pressure on the whole space $ X $.
	
	We denote by $ \mathcal{M}(X) $ the set of all Borel probability measures on $ X $, let $ \mathcal{M}_f(X) $ be the set of $ f $-invariant Borel probability measures and $ \mathcal{M}^e_f(X) $ be the set of ergodic measures in $ \mathcal{M}_f(X) $.
	\subsubsection{\texorpdfstring{Decompositions}{Decompositions}}
	\begin{mdef}
		Given $ \mathcal{D}\subset X\times\mathbb{N} $, a decomposition ($\mathcal{P},\mathcal{G},\mathcal{S}$) for $ \mathcal{D} $ consists of $ \mathcal{P},\mathcal{G},\mathcal{S}\subset X\times\mathbb{N}\cup\{0\} $ and three functions $ p,g,s:\mathcal{D}\to \mathbb{N}\cup\{0\} $ such that for every $ (x,n)\in\mathcal{D} $, the values $ p=p(x,n)$, $g=g(x,n)$, $s=s(x,n) $ satisfy $ n=p+g+s $,  and
		$$
		(x,p)\in \mathcal{P},\quad (f^px,g)\in \mathcal{G}, \quad (f^{p+g}x,s)\in \mathcal{S}.
		$$
	\end{mdef}
	If $ \mathcal{D}=X\times\mathbb{N} $, we call $ (\mathcal{P},\mathcal{G},\mathcal{S}) $ is a decomposition for $ (X,f) $. Given a decomposition $ (\mathcal{P},\mathcal{G},\mathcal{S}) $ and $ N\in \mathbb{N} $, write $ \mathcal{G}^N $ for the set of orbit segments $ (x,n)\in \mathcal{D} $ for which $ p\leq N, s\leq N $.
	We assume that $ (x,0) $ is the empty set. This allows orbit segments to be decomposed in a trivial way, such as choosing a decomposition $ (\mathcal{P},\mathcal{G}) $ where $ \mathcal{S}=\emptyset $.
	\subsubsection{\texorpdfstring{Specification property}{Specification property}}
	\begin{mdef}
		For $ \delta>0 $, we say $ \mathcal{G}\subset X\times \mathbb{N} $ has \textit{$ (W) $-specification at scale $ \delta $} if there exists $ s\in\mathbb{N} $ such that for every $ \{(x_i,n_i)\}_{i=1}^k\subset \mathcal{G} $ there exists a point $ y $ and a sequence of ``gluing times" $ s_1,\cdots,s_k\in \mathbb{N} $ with $ s_i\leq s $ such that writing $ N_j=\sum_{i=1}^j n_i+\sum_{i=1}^{j-1}s_i $, and $ N_0=s_0=0 $, we have
		$$
		d_{n_j}(f^{N_{j-1}+s_{j-1}}y,x_j)<\delta \text{ for every }1\leq j\leq k.
		$$
	\end{mdef}
	\begin{mdef}
		For $ \delta>0 $, we say $ \mathcal{G}\subset X\times \mathbb{N} $ has \textit{specification at scale $ \delta $} if there exists $ p\in\mathbb{N} $ such that for every $ \{(f^{a_i}x_i,b_i-a_i)\}_{i=1}^k\subset \mathcal{G} $ with $ a_1\leq b_1<a_2\leq b_2<\cdots<a_k\leq b_k $ and $ a_i-b_{i-1}\geq p(2\leq i\leq k) $, there exists a point $ y $ such that
		$$
		d(f^{j}y,f^jx_i)<\delta \text{ for every }a_i\leq j\leq b_i, 1\leq i\leq k.
		$$
	\end{mdef}
	It is not difficult to see that the specification property implies the $ (W) $-specification property.
	\begin{mdef}
		Given $ \mathcal{G}\subset X\times\mathbb{N} $ and $ \delta>0 $, define $ \mathcal{G} $ has \textit{tail specification at scale $ \delta $} if there exists $ N_0\in\mathbb{N} $ so that the $ (W) $-specification property holds for orbit segments $ (x_i,n_i)\in\mathcal{G} $ with $ n_i\geq N_0 $.
	\end{mdef}
	
	\subsubsection{\texorpdfstring{The Bowen property}{The Bowen property}}
	\begin{mdef}
		Given $ \varepsilon>0 $, $ \mathcal{C}\subset X\times\mathbb{N} $, we say that $ \phi: X\to \mathbb{R} $ has the \textit{Bowen property } at scale $ \varepsilon $ on $ \mathcal{C} $ if there exists $ K>0 $ such that $ \lvert \Phi_0(x,n)-\Phi_0(y,n)\rvert\leq K $ for all $ (x,n)\in \mathcal{C} $ and $ y\in B_n(x,n) $.
	\end{mdef}
	The constant $ K $ in the above definition is usually called the \textit{distortion constant} for the Bowen property.
	\subsubsection{\texorpdfstring{Obstructions to expansivity}{Obstructions to expansivity}}
	Now $ f $ is a homeomorphism on $ X $. For a fixed $ \varepsilon>0 $ and every $ x\in X $, we can define that
	$$
	\Gamma_{\varepsilon}(x):=\{y\in X\;|\; d_n(x,y)\leq \varepsilon\text{ for all }n\in\mathbb{Z}\}.
	$$
	So if there exists $ \varepsilon>0 $ such that $ \Gamma_{\varepsilon}(x) =\{x\} $ for every $ x\in X $, $ f  $ is expansive at scale $ \varepsilon $. Define the set
	$$
	NE(\varepsilon):=\{x\in X \;|\;\Gamma_{\varepsilon}(x)\neq \{x\}\}
	$$
	and an $ f $-invariant measure $ \mu $ is almost expansive at scale $ \varepsilon $ if $ \mu(NE(\varepsilon))=0 $.
	\begin{mdef}
		For a potential $\phi: X\to \mathbb{R}$ and $\varepsilon>0$, \textit{the pressure of obstructions to expansivity at scale $ \varepsilon $} for a homeomorphism $ f $ is defined as
		\begin{align*}
			P^\perp_{\mathrm{exp}}(f,\phi,\varepsilon)&=\sup\{h_{\mu}(f)+\int \phi d\mu \;|\; \mu\in \mathcal{M}^e_f(X)\text{ and }\mu(NE(\varepsilon))>0\}\\
			&=\sup\{h_{\mu}(f)+\int \phi d\mu \;|\; \mu\in \mathcal{M}^e_f(X)\text{ and }\mu(NE(\varepsilon))=1\}.
		\end{align*}
	\end{mdef}
	Let $ P^\perp_{\mathrm{exp}}(\phi)=\lim\limits_{\varepsilon\to 0}P^\perp_{\mathrm{exp}}(f,\phi,\varepsilon) $. Because $P^\perp_{\mathrm{exp}}(f,\phi,\varepsilon)$ is monotonic when considered as a function of $\varepsilon$ then the limit about $\varepsilon$ exists.
	
	The following results in \cite{CT16} and \cite{Yang} respectively proved the uniqueness of equilibrium states for certain systems.
	\begin{thm} $(\mathrm{Theorem \text{ }5.6 \text{ of } \cite{CT16}})$\label{CT}
		Let $ X $ be a compact metric space, $ f:X\to X $ a homeomorphism and $ \phi:X\to \mathbb{R} $ a continuous potential function. Suppose that there are $ \varepsilon>0$, $\delta>0 $ with $ \varepsilon>40\delta $ such that  $P^\perp_{\mathrm{exp}}(f,\phi,\varepsilon)<P(f,\phi) $ and  there exists $ \mathcal{D}\subset X\times \mathbb{N} $ admits a decomposition $ (\mathcal{P},\mathcal{G},\mathcal{S}) $ with the following properties:
		\begin{flushleft}
			(1) For every $ N\in\mathbb{N} $, $ \mathcal{G}^N  $ has tail $ (W) $-specification at scale $ \delta $;\\
			(2) $ \phi $ has the Bowen property at scale $ \varepsilon $ on $ \mathcal{G} $;\\
			(3) $P(\mathcal{P}\cup\mathcal{S}\cup\mathcal{D}^c,\phi,\delta,\varepsilon)<P(f,\phi) $.\\
			Then there exists a unique equilibrium state for $ (f,\phi) $.
		\end{flushleft}
	\end{thm}
	\begin{thm} $(\mathrm{Theorem \text{ }A} \text{ of }\cite{Yang})$\label{2.2}
		Let $ X $ be a compact metric space, $ f:X\to X $ a homeomorphism and $ \phi:X\to \mathbb{R} $ a continuous potential function. Suppose that there exist $ \varepsilon>0$, $\delta>0 $ with $ \varepsilon\geq 2000\delta $ such that $ P^\perp_{\mathrm{exp}}(f,\phi,\varepsilon)<P(f,\phi) $ and $ X\times \mathbb{N} $ admits a decomposition $ (\mathcal{P},\mathcal{G},\mathcal{S})$ with the following properties:
		\begin{flushleft}
			(1) $ \mathcal{G}  $ has $ (W) $-specification at scale $ \delta $;\\
			(2) $ \phi $ has the Bowen property at scale $ \varepsilon $ on $ \mathcal{G} $;\\
			(3) $ P(\mathcal{P}\cup\mathcal{S},\phi,\delta,\varepsilon)<P(f,\phi) $.\\
			Then there exists a unique equilibrium state for $ (f,\phi) $.
		\end{flushleft}
	\end{thm}
	Theorem \ref{CT} is sometimes referred to as the Climenhaga-Thompson Criterion (abbrev. CT Criterion). As an application of the above criteria, Mongez and Pacifico \cite{Mongez} gave the following result for certain partially hyperbolic diffeomorphisms.
	\begin{thm}$(\mathrm{Theorem \text{ }A} \text{ of }\cite{Mongez})$ \label{Mon}
		Let $f: M\to M$ be a $C^{1+}$ partially hyperbolic diffeomorphism of a compact manifold $M$ with $TM=E^u\oplus E^c\oplus E^s$ and $\phi:M\to \mathbb{R}$ a H\"{o}lder continuous potential. Assume that dim($E^c$)=1 and the unstable foliation $\mathcal{F}^u(f)$ is minimal. If $h^u(f)-h^s(f)>\sup\phi-\inf\phi\geq 0$, then there exists a $C^1$ neighborhood $\mathcal{U}$ of $f$ so that $(g,\phi)$ has a unique equilibrium state for every $C^{1+}$ diffeomorphism $g\in \mathcal{U}$.
	\end{thm}
	The main task of this paper is to obtain the counterpart of Theorem \ref{Mon} in the setting of partially hyperbolic endomorphisms under weaker conditions.

	\subsection{Inverse limit space}
	Now we introduce the inverse limit space for a $ C^2 $ endomorphism $ f $ on the compact  Riemannian manifold $ M $ with metric $ d(\cdot,\cdot) $ induced by the Riemannian metric. The notions following can be found in \cite{Qian}.
	
	Let $ M^f $ denote the subset of $ M^{\mathbb{Z}} $ consisting of all full orbits,
	$$
	M^f:=\{\tilde{x}=(x_i)_{i\in\mathbb{Z}} \;|\; x_i\in M, fx_i=x_{i+1},\forall i\in \mathbb{Z}\},
	$$
	where $ M^{\mathbb{Z}} $ is the infinite product space of $ M $ endowed with the product topology.
	Define metric
	$$
	\tilde{d}(\tilde{x},\tilde{y})=\sum_{n=-\infty}^{\infty}2^{-\lvert n\rvert}d(x_n,y_n)
	$$
	for $ \tilde{x}=(x_n)_{n\in\mathbb{Z}}, \tilde{y}=(y_n)_{n\in\mathbb{Z}}\in M^{\mathbb{Z}}$.
	We call $ M^f $ the \textit{inverse limit space} or \textit{orbit space} of system $ (M,f) $. It is clear that $ M^f $ is a closed subset of $ M^{\mathbb{Z}} $. Let $ \Pi $ denote the natural projection from $ M^f $ to $ M $ such that $ \Pi(\tilde{x}) =x_0$ for every $ \tilde{x}=(x_n)_{n\in\mathbb{Z}}\in M^f $ and $ \tau:M^f\to M^f $ be the left shift map which is a homeomorphism. Note that we have $ \Pi\circ \tau=f\circ \Pi $. Sometimes the system $(M^f,\tau)$ is called the natural extension of $(M,f)$.
	
	Let $ \mathcal{M}_f(M) $ represent the set of $ f $-invariant Borel probability measures on $ M $ and $ \mathcal{M}^e_f(M) $ be the set of ergodic measures in $ \mathcal{M}_f(M) $, and for convenience, we will denote them as $ \mathcal{M}_f $, $ \mathcal{M}^e_f $ in the following, respectively. Similarly, we can define $ \mathcal{M}_{\tau} $ and $ \mathcal{M}^e_{\tau} $. We know that $ \Pi $ induces a continuous map from $ \mathcal{M}_\tau $ to $ \mathcal{M}_f $, usually still denoted by $ \Pi $, that is, for any $ \tau $-invariant Borel probability measure $ \tilde{\mu} $ on $ M^f $, $ \Pi $ maps it to an $ f $-invariant Borel probability measure $ \Pi\tilde{\mu} $ on $ M $ defined by
	$
	\Pi\tilde{\mu}(\psi)=\tilde{\mu}(\psi\circ \Pi),
	$ for every continuous function  $ \psi $  on $  M $.
	
	The following proposition guarantees that $ \Pi $ is a bijection between  $ \mathcal{M}_\tau $  and $ \mathcal{M}_f $.
	\begin{prop}$\mathrm{(Proposition \text{ }I.3.1 \text{ of }\cite{Qian})}$
		Let $ f $ be a continuous map on a compact metric space $ X $. For any $   f $-invariant Borel probability measure $ \mu $ on $ X $, there exists a unique  $ \tau$- invariant Borel probability measure $ \tilde{\mu} $ on $ X^f $ such that $ \Pi\tilde{\mu}=\mu $.
	\end{prop}
	
	In addition, $ \mu $ is ergodic with respect to $ f $ if and only if $ \tilde{\mu} $ is ergodic with respect to $ \tau $.
	Besides, the entropy of the corresponding measure is the same, i.e. $$ h_{\mu}(f)=h_{\tilde{\mu}}(\tau) $$
	whenever $ \Pi\tilde{\mu}=\mu $. For more details, we can refer to Proposition I.3.4 in \cite{Qian}.
	
	Next, we provide a useful fact for the inverse limit space. By the  compactness and $ f $ is continuous on $ M $, it is easily to check that for any $ \varepsilon>0 $,   there are $ \delta>0 $ and an integer $ J $ such that if $ \tilde{x}=(x_n)_{n\in\mathbb{Z}},\tilde{y}=(y_n)_{n\in\mathbb{Z}}\in M^f$ with $ d(x_0,y_0)<\delta $, then $ \tilde{d}(\tau^J\tilde{x},\tau^J\tilde{y})<\varepsilon$.
	\subsection{Partially hyperbolic endomorphisms}
	Let $ f $, $ M $ be as in \thmref{Main}, we give some notions and properties about the partial hyperbolicity and non-degeneration.
	
	Let $ E=\Pi^*TM $ for the pull back bundle of the tangent bundle $ TM $ by the projection $ \Pi:M^f\to M $, the tangent map $ Df $ induces a fiber preserving map on $ E $ with respect to the left shift map $ \tau $, defined by $ \Pi^*\circ Df\circ \Pi_* $ which is still denoted as $ Df $.
	\begin{mdef}
		We say $ f:M\to M $ is a partially hyperbolic endomorphism if there exist a continuous splitting of the pull back bundle $ E $ into three subbundles, that is, $ E(\tilde{x})=E^s(\tilde{x})\oplus E^c(\tilde{x})\oplus E^u(\tilde{x}) $ for all $ \tilde{x}\in M^f $ and constants $ \lambda_s,\lambda_1, \lambda_2,\lambda_u$ and $ C>0 $ with $ 0<\lambda_s<1<\lambda_u$, $\lambda_s<\lambda_1\leq\lambda_2<\lambda_u $ such that for every $ \tilde{x}\in M^f $, we have\\
		(1) $ D_{\tilde{x}}f(E^t(\tilde{x}))=E^t(\tau\tilde{x}) ,t\in\{s,c,u\}$;\\
		(2) $ \lVert D_{\tilde{x}}f^nv^s\rVert\leq C\lambda_s^n\lVert v^s\rVert $, for $ v^s\in E^s(\tilde{x}) $, $ n\geq0 $;\\
		(3) $ C^{-1} \lambda_1^n\lVert v^c\rVert \leq \lVert D_{\tilde{x}}f^nv^c\rVert\leq C\lambda_2^n\lVert v^c\rVert $, for $ v^c\in E^c(\tilde{x}) $, $ n\geq0 $;\\
		(4) $ \lVert D_{\tilde{x}}f^nv^u\rVert\geq C^{-1}\lambda_u^n\lVert v^u\rVert $, for $ v^u\in E^u(\tilde{x}) $, $ n\geq0 $.
	\end{mdef}
	\begin{rmk}
		From the definition above we can see that one of the differences between partially hyperbolic diffeomorphisms and partially hyperbolic endomorphisms is that the unstable and center directions at point $ x $ for endomorphisms depending on the orbit of $ x $ although $ E^s(\tilde{x}) $ is uniquely defined, and the number of $ E^u $ or $ E^c $ at point $ x $ may be more than one. There is an example of an Anosov endomorphism which has infinitely many unstable directions for some orbit constructed by Przytycki; more details can be seen \cite{Przytycki}.
	\end{rmk}
	
	From \cite{Qian, Shu}, we have the following properties. For each $ \tilde{x}=(x_n)_{n\in\mathbb{Z}}\in M^f $ and $ \varepsilon>0 $ sufficiently small, define
	\begin{align*}
		W^u_{\varepsilon}(\tilde{x}):=\{z_0\in M \; |\; &\text{there exists }\tilde{z}\in M^f \text{ with }\Pi\tilde{z}=z_0,\\
		&d(z_{-n},x_{-n})<\varepsilon \text{ for }n\in\mathbb{Z} \text{ and }\limsup_{n\to \infty}\frac{1}{n}\log d(z_{-n},x_{-n})\leq -\log\lambda_u\}
	\end{align*}
	which is said to be a \textit{local unstable manifold} of $ f $ at $ \tilde{x} $. Meanwhile, we can define
	\begin{align*}
		\tilde{W}^u_{\varepsilon}(\tilde{x}):=\{\tilde{z}\in M^f \; |\; \Pi\tilde{z}=z_0\in W^u_{\varepsilon}(\tilde{x}),
		d(z_{-n},x_{-n})<\lambda_u^nd(z_0,x_0) \}.
	\end{align*}
	By Theorem IV.2.2 in \cite{Qian} we get $ \Pi_{\tilde{W}^u_{loc}(\tilde{x})} :\tilde{W}^u_{loc}(\tilde{x}) \to W^u_{loc}(\tilde{x})$ is a bijection.
	We define
	\begin{align*}
		W^u(\tilde{x}):=\{z_0\in M \; |\; &\text{there exists }\tilde{z}\in M^f \text{ with }\Pi\tilde{z}=z_0,\\
		& \text{ and }\limsup_{n\to \infty}\frac{1}{n}\log d(z_{-n},x_{-n})\leq -\log\lambda_u\},
	\end{align*}
	and
	\begin{align*}
		\tilde{W}^u(\tilde{x}):=\{\tilde{z}\in M^f \; |\; \Pi\tilde{z}=z_0\in W^u(\tilde{x}) \text{ with }
		\limsup_{n\to \infty}\frac{1}{n}\log d(z_{-n},x_{-n})\leq -\log\lambda_u \}.
	\end{align*}
	From \cite{Shu} there exists a sequence of $ C^1 $ embedded disks $ \{W_{-n}(\tilde{x})\}_{n=0}^{+\infty} $ in $ M $ such that $ fW_{-n}(\tilde{x})\supset W_{-(n-1)}(\tilde{x}) $ for $ n\in\mathbb{Z}^+ $ and
	\begin{align*}
		W^u(\tilde{x})=\bigcup_{n=0}^{+\infty}f^nW_{-n}(\tilde{x}),
	\end{align*}
	which provides that $ W^u(\tilde{x}) $ is an immersed submanifold of $ M $ tangent at $ \Pi(\tilde{x}) $ to $ E^u(\Pi(\tilde{x})) $. As in \cite{Wang}, define the set $\{W^u(\tilde{x}) \;|\; \tilde{x}\in M^f\}$ by $\mathcal{W}^u$ which is called ${W}^u$-foliation and the set $\{\tilde{W}^u(\tilde{x})\;|\; \tilde{x}\in M^f\}$ by $\tilde{\mathcal{W}}^u$ called the unstable foliation.
	
	The condition in our \thmref{Main} that $ \tilde{W}^u(\tilde{x}) $ is dense in $ M $ for every $ \tilde{x}\in M^f $ is achievable. We give an example below.
	\begin{exa}
		Let $ A=\begin{pmatrix}
			3 & 1\\ 1 & 1
		\end{pmatrix}:\mathbb{T}^2 \to \mathbb{T}^2$ an Anosov endomorphism and $ T:S^1\to S^1 $ the irrational rotation. Then every unstable manifold of $ A\times T $ is dense in $ \mathbb{T}^3 $.
	\end{exa}
	
	When $ f $ is non-degenerate, it is convenient for us to discuss the preimages of $ f $. The reader can refer to \cite{Zhu} for details.
	Define $ f $ has \textit{uniform separation of preimages} if for some $ \varepsilon_0>0 $, $ d(x,y)\leq \varepsilon_0 $ and $ fx=fy $ implies $ x=y $, where $\varepsilon_0$ is said to be an \textit{exponent} of separation for $ f $.
	\begin{lem}$\mathrm{(Lemma \textbf{ }2.6 \text{ of }\cite{Zhu})}$\label{separation}
		If $ f:M\to M$ is a non-degenerate endomorphism, or, more generally, if $ f:M\to M $ is a local homeomorphism, then $ f $ has uniform separation of preimages.
	\end{lem}
	
	\subsection{(Un)Stable entropy and (un)stable pressure}
	In this subsection, we introduce some notions of (un)stable entropy and (un)stable pressure for $ \phi\in  C(M) $ for non-degenerate $C^1$ partially hyperbolic endomorphisms. Let $ \phi^*:=\phi\circ \Pi$ on $ M^f $.
	\subsubsection{\texorpdfstring{Unstable entropy and unstable pressure}{Unstable entropy and unstable pressure}}
	Here we give the definitions of unstable entropy and unstable pressure of $C^1$ partially hyperbolic endomorphisms defined as \cite{Wang}.
	Suppose that $ \eta $ is  a measurable partition of $ M^f $, $ \eta(\tilde{x}) $ means the element in $ \eta $ containing $ \tilde{x} $. Let $ \tilde{\mu} $ be a $ \tau $-invariant Borel probability measure on $ M^f $ with $ \Pi\tilde{\mu}=\mu $.
	\begin{mdef}
		We say a measurable partition $ \eta $ of $ M^f $ is \textit{subordinate to }$ {W}^u $-foliation if for $ \tilde{\mu} $-$a.e.$ $ \tilde{x} $, $ \eta(\tilde{x}) $ has the following properties:\\
		(1) $ \Pi|_{\eta(\tilde{x})}:\eta(\tilde{x})\to \Pi(\eta(\tilde{x})) $ is bijective;\\
		(2) There exists a $ k(\tilde{x}) $-dimensional $ C^1 $ embedded submanifold $ W_{\tilde{x}} $ of $ M  $ with $ W_{\tilde{x}}\subset W^u(\tilde{x})  $ such that $\Pi(\eta(\tilde{x}))\subset W_{\tilde{x}} $ and $ \Pi(\eta(\tilde{x})) $ contains an open neighborhood of $ x_0 $ in $ W_{\tilde{x}} $.
	\end{mdef}
	
	For partition $ \eta $ of $ M^f $ there exists a canonical system $ \{\tilde{\mu}^{\eta}_{\tilde{x}}\}_{\tilde{x}\in M^f} $ of conditional measures of $ \tilde{\mu} $ associated with $ \eta $ satisfying\\
	(1) for every measurable set $ \tilde{B}\subset M^f $, $ \tilde{x}\longmapsto \tilde{\mu}^{\eta}_{\tilde{x}}(\tilde{B})$ is measurable;\\
	(2) $ \tilde{\mu}(\tilde{B})=\int_{M^f} \tilde{\mu}^{\eta}_{\tilde{x}}(\tilde{B})d\tilde{\mu}(\tilde{x}) $.
	
	Take $ \varepsilon>0 $, let $ \alpha $ be a finite partition of $ M^f $ with diam$ (\alpha ):=\sup_{A\in \alpha}$diam$(\Pi(A))\leq \varepsilon $ and denote
	$$
	\eta(\tilde{x})=\alpha(\tilde{x})\cap \tilde{W}^u_{\varepsilon}(\tilde{x})
	$$
	which is a new partition of $ M^f $ finer than $ \alpha $. We denote $ \partial(\Pi(\alpha))=\bigcup_{A\in\alpha}\partial(\Pi(A)) $, then if $ \mu(\partial(\Pi(\alpha)))=0 $, $ \eta $ is a measurable partition subordinate to $ {W}^u$-foliation. Let $ \mathcal{P}^u(M^f) $ be the set of partitions of $ M^f $ subordinate to  ${W}^u$-foliation, which are induced by finite partitions with small diameter.
	Recall that the \textit{information function} of $ \alpha $ with respect to $ \tilde{\mu} $ is given by
	$$
	I_{\tilde{\mu}}(\alpha)(\tilde{x}):=-\log\tilde{\mu}(\alpha(\tilde{x})),
	$$
	and the \textit{entropy} of $ \alpha $ with respect to $ \tilde{\mu} $ is defined by
	$$
	H_{\tilde{\mu}}(\alpha):=\int_{M^f} I_{\tilde{\mu}}(\alpha)(\tilde{x}) d\tilde{\mu}(\tilde{x})=-\int_{M^f}   \log\tilde{\mu}(\alpha(\tilde{x}))d\tilde{\mu}(\tilde{x}).
	$$
	The \textit{conditional information function} of $ \alpha $ with respect to $ \eta $ is defined as
	$$
	I_{\tilde{\mu}}(\alpha|\eta)(\tilde{x}):=-\log\tilde{\mu}^\eta_{\tilde{x}}(\alpha(\tilde{x})),
	$$
	then the \textit{conditional entropy} of $ \alpha $ with respect to $ \eta $ is defined as
	$$
	H_{\tilde{\mu}}(\alpha|\eta):=\int_{M^f} I_{\tilde{\mu}}(\alpha|\eta)(\tilde{x}) d\tilde{\mu}(\tilde{x})=-\int_{M^f}   \log\tilde{\mu}^\eta_{\tilde{x}}(\alpha(\tilde{x}))d\tilde{\mu}(\tilde{x}).
	$$
	Next, we introduce the unstable metric entropy presented in \cite{Wang}.
	\begin{mdef}
		The \textit{conditional entropy} of $ f $ for a finite measurable partition $ \alpha $ of $ M^f $ with respect to $ \eta\in \mathcal{P}^u(M^f) $ is given by
		\begin{align*}
			h_{\mu}(f,\alpha|\eta)=\limsup_{n\to \infty}\frac{1}{n}H_{\tilde{\mu}}(\alpha_0^{n-1}|\eta).
		\end{align*}
		The \textit{conditional entropy} of $ f $  with respect to $ \eta\in \mathcal{P}^u(M^f) $ is given by
		\begin{align*}
			h_{\mu}(f|\eta)=\sup_{\alpha\in \mathcal{P}(M^f)}h_{\mu}(f,\alpha|\eta),
		\end{align*}
		and the \textit{conditional entropy} of $ f $ along $ W^u $-foliation is given by
		\begin{align*}
			h_{\mu}^u(f)=\sup_{\eta\in \mathcal{P}^u(M^f)}h_{\mu}(f|\eta).
		\end{align*}
	\end{mdef}
	The conditional entropy of $ f $ along $ {W}^u $-foliation here can be regarded as the unstable metric entropy for $ f $.
	
	Next, we introduce the definitions of the unstable topological entropy and unstable pressure for a potential function $\phi\in C(M)$ for $f$ using the $ W^u$-separated sets.
	For fixed $ \delta>0 $ and $ \tilde{x}\in M^f $, we denote $ \overline{W^u(\tilde{x},\delta)} $ the $ \delta $-neighborhood of $ x_0 $ in $ W^u(\tilde{x}) $. Let $ d^u_{\tilde{x}} $ be the metric on $ W^u(\tilde{x}) $ induced by the Riemannian structure and set $ d^u_n(y,z):=\max_{0\leq j\leq n-1}\{d^u_{\tau^j\tilde{x}}(f^jy,f^jz)\} $ for any $ y,z\in W^u(\tilde{x}) $. A set $ E\subset\overline{W^u(\tilde{x},\delta)}  $ is said to be $ (n,\varepsilon) $ $ W^u  $-separated if $ \forall y,z\in E$ with $y\neq z $, then $ d^u_n(y,z)>\varepsilon $.
	Denote
	\begin{align*}
		P^u(f,\phi, \tilde{x}, \delta, n,\varepsilon):= \sup&\{\sum_{y\in E} e^{S_n\phi(y)}\;|\;
		E \text{ is an } \\
		&(n,\varepsilon)\text{ }
		W^u\text{-separated
			set of }
		\overline{W^u(\tilde{x},\delta)}\},
	\end{align*}
	where $S_n\phi(y)=\sum_{j=0}^{n-1}\phi(f^j(y))$.
	Denote
	\begin{align*}
		\overline{\tilde{W}^u(\tilde{x},\delta)}=\{\tilde{y}\in M^f \;|\; \Pi(\tilde{y})=y_0\in \overline{W^u(\tilde{x},\delta)} \text{ and } d(y_{-n},x_{-n})\leq \lambda_u^nd(x_0,y_0)\}.
	\end{align*}
	Similarly,  define a subset $ \tilde{E}\subset\overline{\tilde{W}^u(\tilde{x},\delta)}  $ is $ (n,\varepsilon) $ $ W^u  $-separated if $ \forall \tilde{y},\tilde{z}\in \tilde{E}$ with $\tilde{y}\neq\tilde{z} $, then $ \tilde{d}^u_n(\tilde{y},\tilde{z})>\varepsilon $, where $\tilde{d}^u_n(\tilde{y},\tilde{z}):=\max_{0\leq j\leq n-1}\{d^u_{\tau^j\tilde{x}}(\Pi(\tau^j\tilde{y}),\Pi(\tau^j\tilde{z}))\} $. Denote denote
	\begin{align*}
		\tilde{P}(\tau, \phi^*, \tilde{x},\delta, n,\varepsilon):=\sup&\{\sum_{\tilde{y}\in \tilde{E}} e^{S_n\phi^*(\tilde{y})}\;|\; \tilde{E} \text{ is an }\\
		&(n,\varepsilon)\text{ }
		W^u\text{-separated
			set of }\overline{\tilde{W}^u(\tilde{x},\delta)}\},
	\end{align*}
	where $S_n\phi^*(\tilde{y})=\sum_{j=0}^{n-1}\phi^*(\tau^j\tilde{y})$. It is easy to check that $\tilde{P}(\tau, \phi^*, \tilde{x},\delta, n,\varepsilon)=P^u(f,\phi, \tilde{x}, \delta, n,\varepsilon)$.
	\begin{mdef}
		The \textit{unstable pressure} for $ f $ is defined as
		\begin{align*}
			P^u(f,\phi)&=\lim\limits_{\delta\to 0}\sup_{\tilde{x}\in M^f} \lim\limits_{\varepsilon\to 0}\limsup_{n\to \infty}\frac{1}{n}\log P^u(f,\phi, \tilde{x}, \delta, n,\varepsilon)\\
			&=\lim\limits_{\delta\to 0}\sup_{\tilde{x}\in M^f} \lim\limits_{\varepsilon\to 0}\limsup_{n\to \infty}\frac{1}{n}\log \tilde{P}(\tau, \phi^*, \tilde{x},\delta, n,\varepsilon),
		\end{align*}
		and the \textit{unstable topological entropy} of $f$ is defined by
		\begin{align*}
			h^u_{\mathrm{top}}(f)=P^u(f,0).
		\end{align*}
	\end{mdef}
	The relation between the unstable metric entropy and the unstable pressure is given by
	$$
	P^u(f,\phi)=\sup\{h^u_{\mu}(f)+\int\phi d\mu \;|\; \mu\in \mathcal{M}^e_f\}.
	$$
	and in particular, there is a variational principle for unstable topological entropy
	$$
	h^u_{\mathrm{top}}(f)=\sup\{h^u_{\mu}(f) \;|\; \mu\in \mathcal{M}^e_f\}.
	$$

	\subsubsection{\texorpdfstring{Stable entropy and stable pressure}{Stable entropy and stable pressure}}
	Now we give some notions about the stable entropy and stable pressure for $\phi\in C(M)$ of non-degenerate $C^1$ partially hyperbolic endomorphisms in \cite{Zhu, Li}. Similar to the definition of $\mathcal{P}^u$, we can define $ \mathcal{P}^s $ to denote the set of measurable partitions of $ M $ subordinate to the $ {W}^s $-foliation, which are induced by finite partitions with small enough diameter.
	\begin{mdef}
		The \textit{conditional entropy} of $ f $ for a finite measurable partition $ \alpha $ of $ M $ with respect to $ \eta\in \mathcal{P}^s $ is given by
		\begin{align*}
			h_{\mu}(f,\alpha|\eta)=\limsup_{n\to \infty}\frac{1}{n}H_{\mu}(\alpha_0^{n-1}|f^{-n}\eta).
		\end{align*}
		The \textit{conditional entropy} of $ f $  with respect to $ \eta $ is given by
		\begin{align*}
			h_{\mu}(f|\eta)=\sup_{\alpha\in \mathcal{P}}h_{\mu}(f,\alpha|\eta),
		\end{align*}
		and the \textit{stable metric entropy} of $ f $  is given by
		\begin{align*}
			h_{\mu}^s(f)=\sup_{\eta\in \mathcal{P}^s}h_{\mu}(f|\eta).
		\end{align*}
	\end{mdef}
	There are two types of stable topological entropy and stable pressure of $ f $ on $ M $ as follows. Let $ W^s(x,\delta) $ be the open ball inside the stable manifold of $ x $ of radius $ \delta>0$ with respect to the metric $d^s$.
	\begin{mdef}
		The \textit{stable topological pressure} for $ f $ is defined as
		\begin{align*}
			P^s_p(f,\phi)=\lim\limits_{\delta\to 0}\sup_{x\in M}\lim\limits_{\varepsilon\to 0}\limsup_{n\to \infty} \frac{1}{n}\log  P^s(f,\phi,x,\delta,n,\varepsilon)
		\end{align*}
		and
		\begin{align*}
			P^s_m(f,\phi)=\lim\limits_{\delta\to 0}\lim\limits_{\varepsilon\to 0}\limsup_{n\to \infty} \frac{1}{n}\log \sup_{x\in M} P^s(f,\phi,x,\delta,n,\varepsilon)
		\end{align*}
		where $P^s(f,\phi,x,\delta,n,\varepsilon):=\sup_E\{\sum_{y\in E}e^{S_n\phi(y)}\}$ and $E$ ranges over all $(n,\varepsilon)$ $W^s$-separated sets of $f^{-n}\overline{W^s(x,\delta)}$.
		In particular, for the case of $\phi=0$, we obtain the following definitions of the \textit{stable topological entropy} for $f$, i.e.
		\begin{align*}
			h^s_{p,\mathrm{top}}(f):=\lim_{\delta\to 0}\sup_{x\in M}\lim_{\varepsilon\to 0}\limsup_{n\to \infty}\frac{1}{n}\log s(n,\varepsilon,f^{-n}\overline{W^s(x,\delta)})
		\end{align*}
		and
		\begin{align*}
			h^s_{m,\mathrm{top}}(f):=\lim_{\delta\to 0}\lim_{\varepsilon\to 0}\limsup_{n\to \infty}\sup_{x\in M}\frac{1}{n}\log s(n,\varepsilon,f^{-n}\overline{W^s(x,\delta)}),
		\end{align*}
		where $ s(n,\varepsilon,f^{-n}\overline{W^s(x,\delta)}) $ is the maximal cardinality of any $ (n,\varepsilon) $ $ W^s $-separated subset of $ f^{-n}\overline{W^s(x,\delta)} $.
	\end{mdef}
	
	We can also define stable topological pressure using $(n,\varepsilon)$ $W^s$-spanning sets or open covers and more details can be found in \cite{Li}.
	
	There is a variational principle relating stable metric entropy and stable topological entropy, that is,
	$$
	h^s_{p,\mathrm{top}}(f)=h^s_{m,\mathrm{top}}(f)=\sup\{h^s_{\mu}(f) \;|\; \mu\in \mathcal{M}_f\}=\sup\{h^s_{\mu}(f) \;|\; \mu\in \mathcal{M}_f^e\}.
	$$
	Similarly, we have the variational principle about the stable pressure,
	$$
	P^s_p(f,\phi)=P^s_m(f,\phi)=\sup\{h^s_{\mu}(f)+\int\phi d\mu \;|\; \mu\in \mathcal{M}_f^e\}.
	$$
	Let's denote the unstable topological entropy, the stable topological entropy and the stable topological pressure of $ f $ by $ h^u(f) $, $ h^s(f) $ and $P^s(f,\phi)$, respectively, for the sake of argument.
	
	\section{Robustness of the conditions in Main Theorem}
	In this section, we talk about the robustness of the conditions for the non-degenerate $ C^2 $  partially hyperbolic endomorphism $ f $ on $M$ and potential $\phi$ as in  \thmref{Main} and $ \phi^*:=\phi\circ \Pi $. Denote  $ \mathrm{End}^2(M) $ the set of all $ C^2 $ endomorphisms on $ M $ and $ \mathrm{NDPHE}^2(M) $ the set of non-degenerate partially hyperbolic endomorphisms in $ \mathrm{End}^2(M) $.
	
	\subsection{\texorpdfstring{The relationship between unstable pressure and stable pressure}{The relationship between unstable pressure and stable pressure}}
	\begin{prop} \label{3.1}
		If  the potential function $ \phi:M\to \mathbb{R} $ satisfies $ P^u(f,\phi)>P^s(f,\phi) $, then for any $ 0<a<\frac{P^u(f,\phi)-P^s(f,\phi)}{2} $
		there exists a $ C^1 $ neighborhood  $ \mathcal{U} $ of $ f $ such that
		$ P^u(g,\psi)>P^s(g,\psi)$
		for any $\psi\in U_{a}(\phi)$ and any $ g\in \mathrm{NDPHE}^2(M)\cap \mathcal{U} $.
	\end{prop}
	To prove this proposition, we need to prove $ P^u(\cdot,\phi): \mathrm{End}^2(M)\to \mathbb{R} $ is lower semi-continuous and $ P^s(\cdot,\phi): \mathrm{NDPHE}^2(M)\to \mathbb{R} $ is upper semi-continuous at $ f $  with the $ C^1 $ topology, respectively.
	
	According to the result in Theorem $ A$ 1 of \cite{Cantarino}, there exist equilibrium states for any $C^1$ partially hyperbolic endomorphism with one-dimensional center bundle and any continuous potential function.
	
	\begin{lem}
		Let $ P^u(f,\phi)>P^s(f,\phi) $. Then for any ergodic equilibrium state $ \mu $  for $ (f,\phi) $, we have $ h_{\mu}(f)>0 $ and the central Lyapunov exponent  $ \lambda^c(\mu,f)<0 $.
	\end{lem}
	\begin{pf}
		If not, assume there is an ergodic  equilibrium state $ \mu $ for $ (f,\phi) $ with $ h_{\mu}(f)=0 $ or $ \lambda^c(\mu,f)\geq 0$. Since the stable metric entropy $ h^s_\mu(f)\geq 0 $, then
		\begin{align*}
			P(f,\phi)=\int \phi d\mu \leq P^s(f,\phi)<P^u(f,\phi)\leq P(f,\phi),
		\end{align*}
		which is a contradiction.
		
		From Corollary C.1 in \cite{Zhu} we have $ h_{\mu}(f)=h^s_{\mu}(f) $ whenever $ \lambda^c(\mu,f) \geq 0$, then
		\begin{align*}
			P(f,\phi)&=h_{\mu}(f)+\int\phi d\mu=h^s_{\mu}(f)+\int\phi d\mu\\
			&\leq P^s(f,\phi)<P^u(f,\phi)\leq P(f,\phi),
		\end{align*}
		which is also a contradiction.
		This completes the proof of this lemma.
	\end{pf}
	
	Next, we will use the results in \cite{Yun, Chung, Gelfert} that for any $ C^{1+\alpha} $ mapping, there exists a horseshoe whose topological pressure can approximate the sum of the metric entropy and the integral of the potential whenever the metric entropy is positive for an ergodic hyperbolic measure.
	
	\begin{lem}\label{lower}
		The unstable pressure $ P^u(\cdot,\phi): \mathrm{End}^2(M) \to \mathbb{R}$ is lower semi-continuous at $ f  $ with $ C^1 $ topology whenever $ P^u(f,\phi)>P^s(f,\phi) $.
	\end{lem}
	\begin{pf}
		Let $\mu$ be an ergodic equilibrium state for $ (f,\phi) $. By Lemma 3.2, $ \mu $ is an ergodic hyperbolic measure for $ f $ with $ h_{\mu}(f)>0 $. By \cite{Yun, Chung, Gelfert}, for any $ \varepsilon>0 $ there exists a hyperbolic horseshoe $ H_{\varepsilon} $ of $ f $ such that
		\begin{align*}
			P(f|_{H_{\varepsilon}},\phi|_{H_{\varepsilon}})
			>P_{\mu}(f,\phi)-\varepsilon=P^u(f,\phi)-\varepsilon.
		\end{align*}
		
		As we know that the horseshoe for $ C^{1+\alpha} $ mapping is structurally stable under small $ C^1 $ perturbation in \cite{Chung}, i.e. there is a $ C^1 $ neighborhood $ \mathcal{U}_1  $ of $ f $ such that  there exists a hyperbolic horseshoe $ H_g $ and $
		\tau|_{(H_g)^g}: (H_g)^g \to (H_g)^g  $ is topologically conjugate to $\tau|_{(H_\varepsilon)^f}: (H_\varepsilon)^f \to (H_\varepsilon)^f  $ by $ h_g:(H_g)^g \to (H_\varepsilon)^f $ with $ \lVert h_g-id \rVert<\delta $
		for any $ g\in \mathcal{U}_1 $ where $(H_g)^g$, $(H_\varepsilon)^f$ are the inverse limit spaces corresponding to the horseshoes respectively and $ \lVert \phi^*\circ h_g-\phi^*\rVert<\varepsilon $ if $ \lVert h_g-id \rVert<\delta $. Then we have
		\begin{align*}
			P(\tau|_{(H_{\varepsilon})^f},\phi^*|_{(H_{\varepsilon})^f})=P(\tau|_{(H_g)^g}, \phi^*\circ h_g).
		\end{align*}
		and
		\begin{align*}
			P^u(f,\phi)-2\varepsilon<P(g|_{H_g}, \phi|_{H_g}).
		\end{align*}
		In order to relate the above inequality to the unstable pressure of $g$, we need to establish the existence of the equilibrium state for $g|_{H_g}$ as follows.
		\begin{claim}
			If $ H_g $ is a hyperbolic set for the endomorphism $ g $, then $ \tau|_{(H_g)^g} $ is expansive.
		\end{claim}
		In other words, we claim that there is a constant $ \eta>0 $ small enough such that for any $ \tilde{x},\tilde{y} \in (H_g)^g$ with $ \tilde{x}\neq \tilde{y}$ then $ \tilde{d}(\tau^k\tilde{x}, \tau^k\tilde{y})>\eta $ for some $ k\in \mathbb{Z} $. Otherwise $ \tilde{y}\in \tilde{W}^u_{\eta}(\tilde{x})\cap  \tilde{W}^s_{\eta}(\tilde{x})$.
		Since $ \tilde{x}\in  \tilde{W}^u_{\eta}(\tilde{x})\cap  \tilde{W}^s_{\eta}(\tilde{x})$ and $ \Pi_{\tilde{W}^u_{loc}(\tilde{x})}: \tilde{W}^u_{loc}(\tilde{x}) \to W^u_{loc}(\tilde{x}) $  is bijective, then by the local product structure in  Theorem IV.2.3 of \cite{Qian}, we have $\tilde{x}=\tilde{y} $. This completes the proof of this claim.
		
		Since $ H_g $ is hyperbolic, then $\tau|_{(H_g)^g} $ is expansive by Claim 3.4 and there exists an ergodic equilibrium state for $ (\tau|_{(H_g)^g},\phi^*|_{(H_g)^g}) $ denoted by $ \tilde{\mu}_g $, i.e. $ P_{\tilde{\mu}_g}(\tau,\phi^*)=P(g|_{H_g}, \phi|_{H_g})$. 	Note that $ \Pi\tilde{\mu}_g=\mu_g $ is also an ergodic  equilibrium state for $ g|_{H_g} $ with $ \lambda^c(\mu_g,g)<0 $.
		Thus, we have
		\begin{align*}
			P^u(f,\phi)-2\varepsilon<P_{\mu_g}(g,\phi)=P^u_{\mu_g}(g,\phi)\leq P^u(g,\phi)
		\end{align*}
		for any $ g\in \mathrm{End}^2(M)\cap \mathcal{U}_1 $.
		This completes the proof of this lemma.
	\end{pf}
	
	The remaining part in this subsection will extend the results in \cite{Yan21, CYZ} to the non-degenerate partial hyperbolic endomorphisms to prove the upper semi-continuity of stable pressure. In order to transplant the corresponding conclusions in the diffeomorphisms, we will discuss them in the inverse limit spaces.
	
	As we know that both stable and unstable distributions $ E^s, E^u $ are uniquely integrable to the stable and unstable foliations $ \tilde{\mathcal{W}}^s, \tilde{\mathcal{W}}^u $ with $ T\tilde{\mathcal{W}}^s=E^s $ and $ T\tilde{\mathcal{W}}^u=E^u $ respectively for any partially hyperbolic endomorphism. For any $ \tilde{x}\in M^f $ we note that $ E^s(\tilde{x}) $ with $ \Pi(\tilde{x})=x $ only depends on $ x $, however $ E^u(\tilde{x}) $  may depend on the past, then $ \Pi_{\tilde{W}^s_{loc}(\tilde{x})}: {\tilde{W}^s_{loc}(\tilde{x})}\to {{W}^s_{loc}(\tilde{x})}  $ is many-to-one, but the map $ \Pi_{\tilde{W}^u_{loc}(\tilde{x})}:\tilde{W}^u_{loc}(\tilde{x})\to {W}^u_{loc}(\tilde{x}) $ is bijective. Assume the central direction is one-dimensional, then it is integrable but not necessarily uniquely integrable. We denote the local central integrable curve in $ M $ by $ \gamma^c_{loc}(\tilde{x}) $ for any $ \tilde{x}\in M^f $.
	
	Next, we will define the continuity of stable foliations on the $ M $. That is, we talk about the convergence between stable foliations $ \mathcal{W}^s_{n}\to {\mathcal{W}}^s $ when the non-degenerate  $ C^2 $ partially hyperbolic endomorphisms $ f_n\to f $  in the $ C^1 $ topology.
	
	Set $ \sigma$=dim $E^\sigma$ and $ \mathcal{D}^\sigma  $ is the unit disk of   $ \sigma $-dimensional opening in $ \mathbb{R}^\sigma $ where $ \sigma\in \{u,c,s\} $.
	For any $ \tilde{x}\in M^f $ with $ \Pi(\tilde{x})=x $  if there is  a topological embedding $\Phi: \mathcal{D}^{u+1} \times \mathcal{D}^s \to M $ such that $ \Phi(0)=x $, for every  plaque $ P_{z}=\Phi(\{z\}\times \mathcal{D}^s) $ is contained in a certain leaf of $  \mathcal{W}^s $ and
	\begin{align*}
		\Phi(z, \cdot): \mathcal{D}^s \to M, y\mapsto \Phi(z,y)
	\end{align*}
	is a $ C^1 $ embedding and depends on $ z $ continuously in the $ C^1 $ topology, then we denote the image of $ \Phi $ by $ B $  which is called a  \textit{$ {W}^s $-foliation box}  and  the pull back of $ B $ under $ \Pi $ is a \textit{$ \tilde{{W}}^s $-foliation box} which denoted by $ \tilde{B} $.
	Let $ D=\Phi(\mathcal{D}^{u+1}\times \{0\}) $ and $ \tilde{D}\simeq D\times \{\tilde{x}\} $ be the connected component of $ \Pi^{-1}(D) $ containing the point $ \tilde{x} $.
	We also denote $ \{B, \Phi, D\} $ as the foliation box in $ M $ and $ \{\tilde{B}, \Phi, \tilde{D}\} $ as the foliation box in the inverse limit space $ M^f $.

	\begin{rmk}
		If $ D=\Phi(\mathcal{D}^{u+1} \times \{0\})  \subset \bigcup_{w\in \gamma^c_{loc}(\tilde{x})}W^u_{loc}(\tilde{w})$ where the past orbit of $ w $ is  related to $ \tilde{x} $. In fact, $ \Phi $ has determined the backward orbit of the point near  $ x $ since $ W^u, \gamma^c $ may depend on the past. Because $ f $ is non-degenerate, it ensures that the backward orbit of a point within a small neighborhood of $ x $ is quite close to the backward orbit of $ x $ so that a connected component of the pull back $ \tilde{D}\simeq D \times \{\tilde{x}\} $ is also a neighborhood of $ \tilde{x} $.
	\end{rmk}
	
	Since the space $ M $ is compact, there exists a finite open cover consisting of the  $ {W}^s $-foliation boxes $ \{{B}_i, \Phi_i, {D}_i\}_{i=1}^k $. We say  that \textit{$ {\mathcal{W}}_n^s $ converge to $ {\mathcal{W}}^s $} when
	\begin{itemize}
		\item[$\bullet$] for each $ n\in \mathbb{N} $ there is a finite open cover composed of the $ {W}^s_n $-foliation boxes $ \{{B}_i^n, \Phi_i^n, {D}_i \}_{i=1}^k $ of $ M $;
		\item[$ \bullet $] for any $ 1\leq i \leq k $, $ \Phi^n_i: \mathcal{D}^{u+1} \times \mathcal{D}^s \to M $ converge uniformly to $ \Phi_i: \mathcal{D}^{u+1} \times \mathcal{D}^s \to M $ in the $ C^0 $ topology when $ n\to \infty $;
		\item[$ \bullet $] for any $ 1\leq i\leq k $ and any $ z\in \mathcal{D}^{u+1} $, $ \Phi^n_i(z, \cdot): \mathcal{D}^s \to M $ converges to $ \Phi_i(z, \cdot): \mathcal{D}^s \to M $ in the $ C^1 $ topology when $ n\to \infty $.
	\end{itemize}
	
	As we are considering the non-degenerate partially hyperbolic endomorphisms here $ f_n \to f $ in the $ C^1 $ topology, the convergence of the corresponding stable foliations $ \mathcal{W}^s_n \to \mathcal{W}^s $ and  $ \tilde{\mathcal{W}}^s_n \to \tilde{\mathcal{W}}^s $ naturally holds.
	When applying Theorem A in \cite{Yan21} to  $ f^{-1} $, we can get the upper semi-continuity of the stable metric entropy about the $ C^1 $ diffeomorphisms.
	In order to get the corresponding result for an endomorphism that is non-invertible, we construct some partitions that could approximate the stable metric entropy in the inverse limit spaces.
	
	\begin{prop}\label{h^s}
		Assume that $ \{f_n\} $ is a sequence of  non-degenerate $ C^2 $ partially hyperbolic endomorphisms which converge to $ f $ in the $ C^1 $ topology and there are ergodic $ f_n $-invariant measures $ \mu_n $ converge to an ergodic $ f $-invariant measure $ \mu $ in the weak* topology, then
		\begin{align*}
			\limsup_{n\to \infty} h^s_{\mu_n}(f_n)\leq h^s_{\mu}(f).
		\end{align*}
	\end{prop}
	
	Different from the case of diffeomorphisms in \cite{Yan21}, there is also a convergence between the inverse limit spaces with the convergence of a sequence of non-degenerate endomorphisms. By the continuity of the stable foliations we have the foliation boxes  $ \{B^n_i, \Phi^n_i, D_i\}_{i=1}^k $ and $ \{B_i, \Phi_i, D_i\}_{i=1}^k $ in $ M $ and $ \{\tilde{B}^n_i\}_{i=1}^k $ and $ \{\tilde{B}_i\}_{i=1}^k $ in the inverse limit spaces. Take $ \delta=\inf\{\varepsilon_n,\varepsilon_0 \;|\; n\in \mathbb{N}\} $ where $ \varepsilon_n, \varepsilon_0 $ are the exponents of separation for $ f_n $ and $ f  $ respectively. Without losing generality, we assume there are $ \rho>0, \rho'>0 $ small enough such that $2\tilde{\lambda}\rho'<\delta $ and the local product structure property holds in each $ B_i, B^n_i $ whenever the distance between two points is less than $ \rho' $  and every plaque of any foliation box in the space $ M $ has diameter bounded by $ 2\rho' $ and every plaque of any foliation box in the inverse limit space has diameter bounded by $ 2\rho $, where $\tilde{\lambda}=\max_{\tilde{x}\in M^f}\{\lVert D_{\tilde{x}}f^{-1}|_{E^s(x)}\rVert\}$.
	
	For the above $ \rho'>0 $,  we take $ r_0\ll \rho' $ be the Lebesgue number of the open cover $ \{B_i, \Phi_i, {D}_i\}_{i=1}^k  $ of $ M $. Then  we can get a finite partition $ \mathscr{A} $ of $ M $ satisfying diam$(\mathscr{A})  <\min\{\rho', r_0/3, \delta-2\tilde{\lambda}\rho', 1\}$, $\mu(\partial \mathscr{A})=0 $ and $ \mu_n(\partial \mathscr{A})=0 $ for each $ n\in \mathbb{N} $.
	Moreover, the partition $ \mathscr{A} $ induces a function $ \mathscr{I}: M\to  \{1, 2,\cdots, k\} $ such that
	\begin{align*}
		\mathscr{A}(x)\subset  B_{\mathscr{I}(x)},
	\end{align*}
	where $ \mathscr{A}(x) $ is the element of $ \mathscr{A} $ containing $ x $.
	Because the foliations with $ \mathcal{W}^s_n\to \mathcal{W}^s $ as $ n\to \infty $, we have
	\begin{align*}
		\mathscr{A}(x)\subset B^n_{\mathscr{I}(x)}
	\end{align*}
	whenever $ n $ is large enough. After removing finite boxes we could get $ \mathscr{A}(x)\subset B^n_{\mathscr{I}(x)} $ for each $ n\in \mathbb{N} $.
	Then for all $ x\in M $, there exists a unique point $ y\in D_{\mathscr{I}(x)} $ such that $ x\in W^s_{2\rho'}(y) $ and denote
	\begin{align*}
		\mathscr{A}^s(x)=\mathscr{A}(x)\cap W^s_{2\rho'}(y).
	\end{align*}
	Then $ \mathscr{A}^s=\{\mathscr{A}^s(x)\;|\; x\in M\} $ is a measurable partition of $ M $ which is $ \mu $-subordinate to the $ W^s $-foliation. Similarly, we have $ \mathscr{A}^s_n=\{\mathscr{A}^s_n(x)\;|\; x\in M\} $ is a  partition  $ \mu_n $-subordinate to the $ W^s_n $-foliation.

	Before constructing a sequence of new partitions whose conditional entropy can approximate the stable metric entropy,  there is a basic but important fact as follows.  Denote $ f=f_0$, $\tilde{\mu}=\tilde{\mu}_0, \tilde{W}^s_{2\rho}=\tilde{W}^s_{0, 2\rho}$, $ \{\tilde{B}_i\}_{i=1}^k=\{\tilde{B}^0_i\}_{i=1}^k $ and $ \tilde{x}=\tilde{x}^0 $ with $\tilde{x}=(x_i)_{i\in \mathbb{Z}}$, $f(x_i)=x_{i+1}$.
	\begin{lem}\label{4.5}
		For any integer $ n\geq 0, m\geq 0 $ and every $\tilde{x}^n\in M^{f_n} $ with $ \Pi(\tilde{x}^n)=x $,
		\begin{align*}
			\vee_{j=0}^m \tau^{-j}\tilde{\mathscr{A}}^s_n(\tilde{x}^n)=\vee_{j=0}^{m-1} \tau^{-j}\tilde{\mathscr{A}}_n(\tilde{x}^n) \cap \tau^{-m}\tilde{\mathscr{A}}^s_n(\tilde{x}^n),
		\end{align*}
		where $ \tilde{\mathscr{A}}_n=\Pi^{-1}\mathscr{A} $ is a finite measurable partition of $ M^{f_n} $ and $\tilde{\mathscr{A}}^s_n=\Pi^{-1}\mathscr{A}^s_n$.
	\end{lem}
	
	\begin{pf}
		It is obvious that $ \vee_{j=0}^m \tau^{-j}\tilde{\mathscr{A}}^s_n(\tilde{x}^n) \subset \vee_{j=0}^{m-1} \tau^{-j}\tilde{\mathscr{A}}_n(\tilde{x}^n) \cap \tau^{-m}\tilde{\mathscr{A}}^s_n(\tilde{x}^n) $ for each $ \tilde{x}^n\in M^{f_n} $ since $ \tilde{\mathscr{A}}^s_n \succ \tilde{\mathscr{A}}_n$.
		Suppose there exists $ \tilde{y}^n\in \vee_{j=0}^{m-1} \tau^{-j}\tilde{\mathscr{A}}_n(\tilde{x}^n) \cap \tau^{-m}\tilde{\mathscr{A}}^s_n(\tilde{x}^n) $ but $ \tilde{y}^n\notin \vee_{j=0}^m \tau^{-j}\tilde{\mathscr{A}}^s_n(\tilde{x}^n) $, that is, $ \tilde{y}^n\notin \vee_{j=0}^{m-1} \tau^{-j}\tilde{\mathscr{A}}^s_n(\tilde{x}^n) $ which means
		\begin{align*}
			\vee_{j=0}^{m-1} \tau^{-j}\tilde{\mathscr{A}}^s_n(\tilde{y}^n) \neq \vee_{j=0}^{m-1} \tau^{-j}\tilde{\mathscr{A}}^s_n(\tilde{x}^n).
		\end{align*}
		Let $ \tilde{\mathscr{A}}^s_n(\tau^{m-j} \tilde{y}^n)=\tilde{\mathscr{A}}^s_n(\tau^{m-j} \tilde{x}^n) $ with $ 0\leq j <k $ and $ \tilde{\mathscr{A}}^s_n(\tau^{m-k} \tilde{y}^n)\neq \tilde{\mathscr{A}}^s_n(\tau^{m-k} \tilde{x}^n) $. Then there exists $z_0\in D_{\mathscr{I}(f^{m-k+1}_nx)}$ such that $ f^{m-k+1}_n(y)\in \mathscr{A}(f^{m-k+1}_n(x)) \cap W^s_{n, 2\rho'}(z_0) $ and the diameter of $ \mathscr{A} $ is small sufficiently, we have
		\begin{align*}
			f^{m-k}_n(y)\in \mathscr{A}(f^{m-k}_n(x))\cap W^s_{n, 2\rho'}(z_{-1}),
		\end{align*}
		where $f^{m-k}_n x\in W^s_{n,2\rho'\tilde{\lambda}}(z_{-1})$ with $f(z_{-1})=z_0$.
		By the definition of $ \tilde{\mathscr{A}}^s_n $, we get
		$ \tau^{m-k} \tilde{y}^n \in \tilde{\mathscr{A}}^s_n(\tau^{m-k} \tilde{x}^n)$ which contradicts the hypothesis.
	\end{pf}
	
	Recall that $ \Pi : M^{f_n}\to M $ is the natural projection map for each $ n\geq 0 $. For notational convenience, we denote both maps by $\Pi$ even though they are defined on different domains, respectively, if there is no confusion.
	Denote  $\tilde{\mu}_n=\mu_n \circ \Pi$.  Then by Proposition 3.8 in \cite{Zhu} and Lemma 3.6 we have
	\begin{align*}
		h^s_{\mu_n}(f_n)=\lim\limits_{m\to\infty} \frac{1}{m}H_{\tilde{\mu}_n}(\vee_{j=0}^{m-1}\tau^{-j}\tilde{\mathscr{A}}^s_n|\tau^{-m}\tilde{\mathscr{A}}^s_n)
		=\inf_{m\geq 1}\frac{1}{m}H_{\tilde{\mu}_n}(\vee_{j=0}^{m-1}\tau^{-j}\tilde{\mathscr{A}}^s_n|\tau^{-m}\tilde{\mathscr{A}}^s_n).
	\end{align*}
	
	For any $ n\geq 0 $ and each $ 1\leq i \leq k $ , let $ \mathscr{C}_{i,1}\preceq \mathscr{C}_{i,2}\preceq \cdots $ be a sequence of finite partitions on $ {D}_i $ such that
	\begin{itemize}
		\item[$ \bullet $] diam$ (\mathscr{C}_{i,t})\to 0 $ as $ t\to \infty $;
		\item[$ \bullet $] $ {\mu}_{n,i}(\partial \mathscr{C}_{i,t})=0 $ for any $ t\in \mathbb{N} $, where $ {\mu}_{n,i} $ is the projection of $ {\mu}_n $  restricted to $ \cup_{{y}\in {D}_i}{W}^s_{n,2\rho}({y}) $ on $ {D}_i $.
	\end{itemize}
	Let $ \hat{D}_i={D}_i\setminus \cup_{t}\partial\mathscr{C}_{i,t} $ for any $ 1\leq i\leq k $. For each $ i\in \{1,\cdots, k\}$, fix $  t\in \mathbb{N} $ and $ C\in \mathscr{C}_{i,t} $, for any $ x\in C\cap \hat{D}_i $ there is a $ r_{{x},t}>0 $ such that the distance between $ {W}^s_{n, 2\rho'}({x}) $ and any other $ {W}^s_{n, 2\rho'}({y}) $ with $ {y}\in \hat{D}_i\setminus C $ is larger than $ 2r_{{x},t} $ for any $ n\geq 0 $.
	Define the set
	\begin{align*}
		P^n_C=\bigcup_{{x}\in C\cap \hat{D}_i} \bigcup_{{y}\in B({x}, r_{{x},t})\cap {D}_i} {W}^s_{n, 2\rho'}({y}),
	\end{align*}
	we note that $ P^n_C\cap {D}_i $ is open in $ C $ and has full $ {\mu}_{n,i} $ measure in $ C $ and $ {\mu}_{n,i}(\partial(P^n_C \cap {D}_i))=0 $. Then there is a finite partition on the  $ {B}^n_i $ induced by $ \mathscr{C}_{i,t} $ which is given by
	\begin{align*}
		\mathscr{P}_{n,i,t}=\{P^n_C\;|\; C\in \mathscr{C}_{i,t}\}\cup \{{B}^n_i\setminus \bigcup_{C\in \mathscr{C}_{i,t}} P^n_C\}.
	\end{align*}
	For each $ x\in \hat{D}_i $, since diam$( \mathscr{C}_{i,t} ) \to 0 $ we can take a monotonically decreasing sequence of $ r_{{x},t} $ satisfying $ r_{{x},t}\to 0 $ as $ t\to \infty $, there is
	$ 	\mathscr{P}_{n,i,t} \prec \mathscr{P}_{n,i,t+1} $
	in the sense of modulo a  zero $ {\mu}_n $-measure  set and $ \cap_{t}\mathscr{P}_{n,i,t}(x)=W^s_{n, 2\rho'}(x) $ for any $ x\in \hat{D}_i $.
	
	Let $ \Pi^{-1}\mathscr{P}_{n,i,t}=\tilde{\mathscr{P}}_{n,i,t}=\{\tilde{P}^n_C \;|\; C\in \mathscr{C}_{i,t}\}\cup \{\tilde{B}^n_i\setminus \bigcup_{C\in \mathscr{C}_{i,t}}\tilde{P}^n_C\} $ be the partition of $ \tilde{B}^n_i $ in the inverse limit space $ M^{f_n} $ for any $ n\geq 0 $ and $ 1\leq i\leq k $.
	Fix $ n\geq 0 $, for  $ \tilde{\mu}_n$-$a.e.$ $ \tilde{x}^n\in M^{f_n} $ with $ \Pi(\tilde{x}^n)=x $ and any $ \tilde{A}_n\in \tilde{\mathscr{A}}_n $ with $ \tilde{A}_n(\tilde{x}^n)\subset \tilde{B}^n_{\mathscr{I}(x)} $, the partition $ \tilde{\mathscr{P}}_{n,\mathscr{I}(x),t} $ induces a finite measurable partition for $ \vee_{j=0}^m \tau^{-j} \tilde{\mathscr{A}}_n(\tilde{x}^n) $ which is denoted by
	\begin{align*}
		\mathscr{Q}_t=\{\vee_{j=0}^m \tau^{-j} \tilde{\mathscr{A}}_n(\tilde{x}^n) \cap \tilde{P} \;|\;  \tilde{P}\in \tilde{\mathscr{P}}_{n,\mathscr{I}(x),t}\},
	\end{align*}
	and for any $ n\geq 0, m\geq 0 $ and $ t\in \mathbb{N}$,
	\begin{align*}
		\tilde{\mathscr{A}}^m_{n,t}=\{\mathscr{Q}_t \;|\; \tilde{x}^n\in M^{f_n}\}
	\end{align*}
	can be regarded as  a partition of $ M^{f_n} $  and has the following properties: for any $ n\geq 0, m\geq 0 $
	\begin{itemize}
		\item[$ \bullet $] $ \vee_{j=0}^m \tau^{-j}\tilde{\mathscr{A}}_n \prec \tilde{\mathscr{A}}^m_{n,t} \prec \vee_{j=0}^m \tau^{-j}\tilde{\mathscr{A}}^s_{n}$;
		\item[$ \bullet $] $\tilde{\mathscr{A}}^m_{n,t}\nearrow \vee_{j=0}^m \tau^{-j}\tilde{\mathscr{A}}^s_{n}$ as $ t\to \infty $;
		\item[$ \bullet $] $ \tilde{\mu}_n(\partial \tilde{\mathscr{A}}^m_{n,t})=0 $;
		\item[$ \bullet $] $ H_{\tilde{\mu}_n}(\tilde{\mathscr{A}}^{m-1}_{n,t} \;|\;\tau^{-m}\tilde{\mathscr{A}}^0_{n,t})\searrow H_{\tilde{\mu}_n}(\vee_{j=0}^{m-1} \tau^{-j}\tilde{\mathscr{A}}^s_n \;|\;\tau^{-m}\tilde{\mathscr{A}}^s_n) $ when $ t \to \infty $,
	\end{itemize}
	the first three items follow from the construction of $ \tilde{\mathscr{A}}^m_{n,t} $, and the proof of the last one is analogous to that in \cite{Yan21} by applying Lemma \ref{4.5}, the reader can refer to the proof of Proposition 4.5 in \cite{Yan21} for more details.
	
	Now let's finish the proof of Proposition \ref{h^s}. For any $ \varepsilon>0 $, we can choose $ m>0 $ large enough such that
	\begin{align}
		h^s_{\mu}(f)> \frac{1}{m}H_{\tilde{\mu}}(\vee_{j=0}^{m-1}\tau^{-j}\tilde{\mathscr{A}}^s \;|\; \tau^{-m}\tilde{\mathscr{A}}^s)-\frac{\varepsilon}{3}.
	\end{align}
	Then there exists  a large integer $ t=t(m)>0 $ such that for $ n=0 $,
	\begin{align}
		H_{\tilde{\mu}}(\vee_{j=0}^{m-1}\tau^{-j}\tilde{\mathscr{A}}^s \;|\; \tau^{-m}\tilde{\mathscr{A}}^s)>H_{\tilde{\mu}}(\tilde{\mathscr{A}}^{m-1}_{0,t} \;|\; \tau^{-m}\tilde{\mathscr{A}}^0_{0,t})-\frac{\varepsilon}{3}.
	\end{align}
	Since $ f_n\to f $ in the $ C^1 $ topology and the stable foliations $ \mathcal{W}^s_n $ converge to $ \mathcal{W}^s $, by the properties of partitions $ \tilde{\mathscr{A}}^{m-1}_{n,t} $ for each $ Q_0\in \tilde{\mathscr{A}}^{m-1}_{0,t} $ there exists a sequence of sets $ \{Q_n\in \tilde{\mathscr{A}}^{m-1}_{n,t} \;|\; n\in \mathbb{N}\} $ such that
	\begin{itemize}
		\item[$ \bullet $] $ \overline{int(Q_n)}= \overline{Q_n} $ for any $ n\geq 0 $;
		\item[$ \bullet $] for every compact set $ K\subset int(Q_0) $ there is $ K\subset int(Q_n) $ for large $ n\in \mathbb{N} $ and replacing $ Q_0$, $Q_n $  with their complements also holds;
		\item[$ \bullet $] $\tilde{\mu}(\partial Q_0)=0$ and $ \tilde{\mu}_n(\partial Q_n)=0 $ for any $ n\in \mathbb{N} $.
	\end{itemize}
	Then $ \overline{Q_n}\to \overline{Q_0} $ in the Hausdorff topology so that $ \Pi(\overline{Q_n})\to \Pi(\overline{Q_0}) $ in the Hausdorff topology and $ \mu_n(\partial \Pi({Q_n}))=0,  \mu(\partial \Pi({Q_0}))=0$, because $ \mu_n\to\mu $ in the weak* topology and $ \tilde{\mu}_n({Q_n})=\mu_n\circ \Pi({Q_n}) $ for each $ n\geq 0 $, we have
	\begin{align*}
		\lim\limits_{n\to\infty}\tilde{\mu}_n({Q_n})=\tilde{\mu}({Q_0}).
	\end{align*}
	Repeat the same operation to take out the corresponding elements in $ \tau^{-m}\tilde{\mathscr{A}}^0_{n,t} $ and $\tau^{-m}\tilde{\mathscr{A}}^0_{0,t} $. Calculate the conditional entropy, and we get
	\begin{align*}
		\lim\limits_{n\to\infty}H_{\tilde{\mu}_n}(\tilde{\mathscr{A}}^{m-1}_{n,t} \;|\; \tau^{-m}\tilde{\mathscr{A}}^0_{n,t})
		=H_{\tilde{\mu}}(\tilde{\mathscr{A}}^{m-1}_{0,t} \;|\; \tau^{-m}\tilde{\mathscr{A}}^0_{0,t}).
	\end{align*}
	Then there exists large $ N\in \mathbb{N} $ such that for any $ n\geq N $,
	\begin{align*}
		H_{\tilde{\mu}_n}(\tilde{\mathscr{A}}^{m-1}_{n,t} \;|\; \tau^{-m}\tilde{\mathscr{A}}^0_{n,t})-\frac{\varepsilon}{3}
		\leq H_{\tilde{\mu}}(\tilde{\mathscr{A}}^{m-1}_{0,t} \;|\; \tau^{-m}\tilde{\mathscr{A}}^0_{0,t}).
	\end{align*}
	By (1) and (2),  for any $ n\geq N $ we have
	\begin{align*}
		h^s_{\mu}(f)>\frac{1}{m}	H_{\tilde{\mu}_n}(\tilde{\mathscr{A}}^{m-1}_{n,t} \;|\; \tau^{-m}\tilde{\mathscr{A}}^0_{n,t})-\varepsilon>h^s_{\mu_n}(f_n) -\varepsilon.
	\end{align*}
	It follows from the arbitrariness of $ \varepsilon $ that we conclude the proof of Proposition \ref{h^s}.
	
	\begin{cor}\label{upper}
		The stable pressure $ P^s(\cdot,\phi): \mathrm{NDPHE}^2(M) \to \mathbb{R} $ is upper semi-continuous at $f$ with $C^1$ topology.
	\end{cor}
	\begin{pf}
		It is not difficult to verify this from the upper semi-continuity of the stable entropy map $ \mu\mapsto h^s_{\mu}(f) $ at $ \mu\in \mathcal{M}^e_f$, the weak* continuity argument and Proposition \ref{h^s}.
	\end{pf}
	
	\begin{rmk}
		By the results about the unstable pressure for partially hyperbolic endomorphisms in \cite{Wang}, we can obtain that the unstable pressure $P^u(\cdot,\phi): \mathrm{NDPHE}^2(M) \to \mathbb{R}$ is also upper semi-continuous at $f$ with $C^1$ topology through a similar discussion.
	\end{rmk}
	
	\begin{pot3.1}
		For any $ 0<a<\frac{P^u(f,\phi)-P^s(f,\phi)}{2} $ and every $ \psi\in U_a(\phi)=\{\psi\;|\; \lVert \psi-\phi\rVert<a\} $.
		On the one hand, from Lemma \ref{lower} we get $ P^u(\cdot,\phi):\mathrm{End}^2(M)\to \mathbb{R} $ is lower semi-continuous at $ f $ with $ C^1 $ topology, then for any $ \varepsilon>0 $, we can choose a $ C^1 $ neighborhood $ \mathcal{U}_1 $ of $ f $ such that for every $ g\in \mathrm{End}^2(M)\cap \mathcal{U}_1 $, we have
		\begin{align*}
			P^u(g,\psi)+a\geq P^u(g,\phi)>P^u(f,\phi)-\varepsilon.
		\end{align*}
		On the other hand, by Corollary \ref{upper} we have $ h^s(\cdot): \mathrm{NDPHE}^2(M)\to \mathbb{R} $ is upper semi-continuous at $ f $ with $ C^1 $ topology, then for the  $ \varepsilon $ given above, there is a $ C^1 $ neighborhood $ \mathcal{U}_2 $ of $ f $ so that for any $ g\in \mathrm{NDPHE}^2(M)\cap \mathcal{U}_2 $, we get
		\begin{align*}
			P^s(g,\psi)-a\leq P^s(g,\phi)<P^s(f,\phi)+\varepsilon.
		\end{align*}
		Fix $ 0<\varepsilon<\frac{P^u(f,\phi)-P^s(f,\phi)}{2}-a $ and $ \mathcal{U}=\mathcal{U}_1\cap \mathcal{U}_2 $, for each $ \psi\in U_a(\phi) $ and any $ g\in \mathrm{NDPHE}^2(M)\cap \mathcal{U} $ there is
		\begin{align*}
			P^u(g,\psi)-P^s(g,\psi)
			\geq
			P^u(f,\phi)-P^s(f,\phi)-2a-2\varepsilon
			>0.
		\end{align*}
		This completes the proof of this proposition.
	\end{pot3.1}
	
	\subsection{\texorpdfstring{The minimality of the unstable foliation}{The minimality of the unstable foliation}}
	
	For a partially hyperbolic endomorphism $ f: M\to M $, we say the unstable foliation $ \tilde{\mathcal{W}}^u $ is \textit{minimal} if $ \tilde{W}^u(\tilde{x}) $ is dense in $ M^f $ for all $ \tilde{x}\in M^f $.
	The unstable foliation $ \tilde{\mathcal{W}}^u $ is \textit{$\varepsilon$-minimal} if there exists $ R>0 $ such that if $ D $ is a disk contained in a certain leaf of $ \tilde{\mathcal{W}}^u $ with an internal radius $ \eta $ larger than $ R $ then $ D $ is $ \varepsilon $-dense in $ M^f $, where $ D(\tilde{x})=\{\tilde{y}\in \tilde{{W}}^u(\tilde{x}) \;|\; \Pi(\tilde{y})=y\in W^u(\tilde{x}) \text{ with } d_{W^u}(x,y)\leq \eta\} $ and $ d_{W^u}(x,y)=\inf\{\text{length}(\gamma) \;|\; \gamma \text{ is a path from } x \text{ to } y \text{ that is wholly contained in } W^u(\tilde{x})\} $ if we assume the disk $ D $ is centered at $ \tilde{x} $ with $ \Pi(\tilde{x})=x $.
	
	From the unstable manifold theorem for partially hyperbolic endomorphisms we know that if non-degenerate map $ f: M\to M $ with the minimal unstable foliation, then for any $\varepsilon>0$ there exists a $ C^1 $ neighborhood $ \mathcal{U} $ of $ f $ such that the unstable foliation of $ g $ is $ \varepsilon$-minimal for each non-degenerate map $ g\in \mathcal{U} $.
	
	It is well known that every leaf of $ \mathcal{W}^u $ is dense in $ M $ whenever each leaf of $ \tilde{\mathcal{W}}^u $ is dense in $ M^f $. There are similar properties for $ \varepsilon$-minimality.

	\section{Proof of Main Theorem}
	
	\subsection{\texorpdfstring{Estimation of the central Lyapunov exponent of the ergodic equilibrium state}{Estimation of the central Lyapunov exponent of the ergodic equilibrium state}}

	Now we assume that $f $ is a $ C^2 $  non-degenerate partially hyperbolic endomorphism on the compact Riemannian manifold with dim$(E^c)=1 $ and $ \phi $ is a H\"{o}lder continuous potential function on $ M $. Recall that let $\phi^*=\phi\circ \Pi$. Given an ergodic measure $ \mu\in \mathcal{M}^e_f $, we can consider the central Lyapunov exponents by Theorem I.2.1 and Proposition I.3.5 of \cite{Qian}. We use $\tilde{\mu} $ to denote the ergodic $ \tau $-invariant Borel probability measure on $ M^f $ such that $ \Pi \tilde{\mu}=\mu $. We denote $ \lambda^c(\tilde{\mu},\tau) $ (resp. $\lambda^c(\mu,f)$) be the central Lyapunov exponent of $ \tilde{\mu} $ (resp. $ \mu$) whenever $ \tilde{\mu}$ (resp. $\mu $) is ergodic, i.e. for $\tilde{\mu}$-$a.e.$ $ \tilde{x}=(x_n)_{n\in\mathbb{Z}}\in M^f$ (resp. $ \mu$-$a.e.$ $x_0\in M $),  we have
	\begin{align*}
		\lambda^c(\mu,f)=\lambda^c(\tilde{\mu},\tau)
		&=\lim\limits_{m\to\infty}\frac{1}{m}\log\lvert D_{x_0}f^m \xi\rvert\\
		&=\lim\limits_{m\to\infty}\frac{1}{m}\log\lVert D_{x_0}f^m |_{E^c(\tilde{x})}\rVert\\
		&=\lim\limits_{m\to\infty}\frac{1}{m}\log\lVert D_{x_{m-1}}f\circ \cdots \circ D_{x_1}f \circ D_{x_0}f |_{E^c(\tilde{x})}\rVert\\
		&=\lim\limits_{m\to\infty}\frac{1}{m}\sum_{k=0}^{m-1}\log\lVert D_{x_k}f |_{E^c(\tau^k\tilde{x})}\rVert\\
		&=\int_{M^f} \log\lVert Df |_{E^c(\tilde{x})}\rVert d\tilde{\mu},		
	\end{align*}
	where  $ \tilde{x}\in M^f $ with $ \Pi(\tilde{x})=x_0 $, $ 0\neq\xi\in E^c(\tilde{x}) $, and define $\tilde{\varphi}^c(\tilde{x}):=\log\lVert Df |_{E^c(\tilde{x})}\rVert$.
	
	By Corollary C.1 of \cite{Zhu}  we have
	
	\begin{equation}\label{4.1}
		h_{\mu}(f)=h_{\mu}^{s}(f)-\lambda^c \gamma^c,
	\end{equation}
	where $ \lambda^c $ is the negative central Lyapunov exponent and $ \gamma^c $ is the transversal Hausdorff dimension of $ \mu $ on  Pesin center manifold.
	
	\begin{lem}\label{lip}
		If $ P^u(f,\phi)>P^s(f,\phi)$, then for any ergodic equilibrium state $ \tilde{\mu} $ for $ (\tau,\phi^*) $, we have $ \lambda^c(\tilde{\mu}, \tau)<0 $.
	\end{lem}
	\begin{pf}
		Assume that there is an ergodic equilibrium state $\tilde{\mu}$ for $(\tau, \phi^*)$ with $\lambda^c(\tilde{\mu},\tau)\geq 0 $. From Proposition I.3.4 in \cite{Qian}, for the ergodic measure $\mu$ on $M$ with $\Pi\tilde{\mu}=\mu$, there is $h_{\tilde{\mu}}(\tau)=h_{\mu}(f)$. In particular, $P(\tau,\phi^*)=P(f,\phi)$.
		Then by (3) we have,
		\begin{align*}
			h_\mu(f)+\int \phi d\mu
			&=h_\mu^s(f)+\int \phi d\mu
			\leq P^s(f,\phi)
			\\&<P^u(f,\phi)
			\leq P(f,\phi),
		\end{align*}
		which contradicts that $ \mu $ is also an equilibrium state for $ (f,\phi) $.
	\end{pf}
	
	\subsection{\texorpdfstring{A decomposition for $(M^f,\tau)$}{A decomposition for $(M^f,\tau)$}}
	In this subsection, we want to find a decomposition $(\tilde{\mathcal{P}}, \tilde{\mathcal{G}}, \tilde{\mathcal{S}})$ which satisfies the conditions in Theorem \ref{2.2}. In order to show that there exists the specification property on $\tilde{\mathcal{G}}$, we will show that the orbit segments in the good collection $ \mathcal{G}\subset M\times \mathbb{N} $ experience uniform contraction in the center-stable direction.
	
	Define
	\begin{align*}
		P^+(\tau,\phi^*)=\sup\{h_{\tilde{\mu}}(\tau)+\int\phi^* d\tilde{\mu}: \tilde{\mu}\in \mathcal{M}_\tau^e, \lambda^c(\tilde{\mu},\tau)\geq 0\},\\
		P^-(\tau,\phi^*)=\sup\{h_{\tilde{\mu}}(\tau)+\int\phi^* d\tilde{\mu}: \tilde{\mu}\in \mathcal{M}_\tau^e, \lambda^c(\tilde{\mu},\tau)\leq 0\},
	\end{align*}
	where $\mathcal{M}_\tau$ is the set of all $\tau$-invariant probability measures on $M^f$ and $\mathcal{M}^e_\tau$ consists of all ergodic measures in $\mathcal{M}_\tau$.
	We get immediately there is no ergodic equilibrium state $ \tilde{\mu} $ with $ \lambda^c(\tilde{\mu},\tau)\geq 0 $ whenever $ P^u(f,\phi)>P^s(f,\phi) $ by Lemma \ref{lip}. Since $ P(\tau,\phi^*)=\max\{ P^+(\tau,\phi^*),P^-(\tau,\phi^*)\} $, then $ P^+(\tau,\phi^*)<P^-(\tau,\phi^*) $ under the condition stated in Lemma \ref{lip}.

	Next, in order to find a good collection of the orbit segments, we will consider the central Lyapunov exponents, which obviously have the contraction property.
	
	\begin{lem}
		There exists $ r>0 $ such that
		$ \sup\{h_{\tilde{\mu}}(\tau)+\int\phi^* d\tilde{\mu}: \tilde{\mu}\in \mathcal{M}_\tau, \lambda^c(\tilde{\mu},\tau)\geq -r\}<P(\tau,\phi^*) $.
	\end{lem}
	\begin{pf}
		First, we claim that
		$ \sup\{h_{\tilde{\mu}}(\tau)+\int\phi^* d\tilde{\mu}: \tilde{\mu}\in \mathcal{M}_\tau, \lambda^c(\tilde{\mu},\tau)\geq 0\}<P(\tau,\phi^*) $. If not, we can assume there exists an invariant probability measure $ \tilde{\mu} $ with $ \lambda^c(\tilde{\mu},\tau)\geq 0 $, but $ \tilde{\mu} $ is an equilibrium state for $ (\tau,\phi^*) $. By the ergodic decomposition theorem, we have
		$ \lambda^c(\tilde{\mu},\tau)=\int \tilde{\varphi}^c d\tilde{\mu}=\int_{\mathcal{M}_\tau^e} (\int \tilde{\varphi}^c dm) d\tau(m) \geq 0$,
		and by  Theorem 8.4 in \cite{Walters} we get
		$ h_{\tilde{\mu}}(\tau)=\int_{\mathcal{M}_\tau^e} h_m(\tau)d\tau(m) $.
		Then
		$$
		P(\tau,\phi^*)=h_{\tilde{\mu}}(\tau)+\int \phi^* d\tilde{\mu}
		=\int_{\mathcal{M}_\tau^e} (h_m(\tau)+\int \phi^* dm)d\tau(m),
		$$
		so each $ m\in \mathcal{M}_\tau^e $ must be the equilibrium state for $ (\tau,\phi^*) $. Therefore, there exists a measure $ m_0\in\mathcal{M}_\tau^e  $ such that
		$ \lambda^c(m_0,\tau)=\int \tilde{\varphi}^c dm_0\geq 0 $
		and
		$ P(\tau,\phi^*)=h_{m_0}(\tau)+\int \phi^* dm_0 $,
		which contradicts with  $ P^+(\tau,\phi^*)<P^-(\tau,\phi^*) $.

		Let's continue to prove Lemma 4.2. Suppose that for any $ r>0 $, we have
		$$
		\sup\{h_{\tilde{\mu}}(\tau)+\int \phi^* d\tilde{\mu}:\tilde{\mu} \in \mathcal{M}_\tau,\lambda^c(\tilde{\mu},\tau)\geq -r\}=P(\tau,\phi^*).
		$$
		Then we can choose a sequence $ \{r_n\}_{n\geq 1} $ with $ \lim\limits_{n\to\infty}r_n=0 $ and
		$$
		\sup\{h_{\tilde{\mu}}(\tau)+\int \phi^* d\tilde{\mu}:\tilde{\mu} \in \mathcal{M}_\tau,\lambda^c(\tilde{\mu},\tau)\geq -r_n\}=P(\tau,\phi^*).
		$$
		For any $ n\in\mathbb{N}  $, there exist $ \tilde{\mu}_n\in\mathcal{M}_\tau $ with $ \lambda^c(\tilde{\mu}_n,\tau)\geq -r_n $ and
		$$
		h_{\tilde{\mu}_n}(\tau)+\int \phi^* d\tilde{\mu}_n\geq P(\tau,\phi^*)-\frac{1}{n}.
		$$
		Since $ \mathcal{M}_\tau $ is compact, let $ \tilde{\mu}' $ be a limit point of $ \tilde{\mu}_n $, we have $ \tilde{\mu}' $ is $ \tau $-invariant and $ \lambda^c(\tilde{\mu}',\tau)\geq 0 $ by the weak* topology. Moreover, $ P(\tau,\phi^*)\leq \limsup\limits_{n\to\infty}(h_{\tilde{\mu}_n}(\tau)+\int \phi^*d\tilde{\mu}_n)\leq h_{\tilde{\mu}'}(\tau)+\int \phi^*d\tilde{\mu}' $
		which contradicts the claim above.
	\end{pf}
	Now we can define a decomposition $ (\tilde{\mathcal{P}}, \tilde{\mathcal{G}}, \tilde{\mathcal{S}}) $ for $ M^f\times\mathbb{N} $ and describe a good collection of orbit segments $ \mathcal{G}\subset M\times\mathbb{N} $. We take $ \tilde{\mathcal{S}}=\emptyset $, fix $ r>0 $ as defined in Lemma 4.2, define
	\begin{align*}
		\tilde{\mathcal{P}}:=\{(\tilde{x},n)\in M^f\times\mathbb{N}\; | \; S_n\tilde{\varphi}^c(\tilde{x})\geq -rn \}		
	\end{align*}
	and
	\begin{align*}
		\tilde{\mathcal{G}}:=\{(\tilde{x},n)\in M^f\times\mathbb{N}\;|\; S_j\tilde{\varphi}^c(\tilde{x})< -rj \text{ for all } 1\leq j \leq n \},
	\end{align*}
	where $ S_j\varphi^c(x)=\sum_{i=0}^{j-1}\tilde{\varphi}^c(\tau^{i}\tilde{x}) $. For any orbit segment $ (\tilde{x},n)\in M^f\times \mathbb{N} $, we set $ p=p(\tilde{x},n) $ be the maximal integer such that $ (\tilde{x},p) $ in the bad collection, that is, $ (\tilde{x},p)\in \tilde{\mathcal{P}} $. Then we can easily check that $ (\tau^p\tilde{x},n-p)\in \tilde{\mathcal{G}} $ and $ (\tilde{\mathcal{P}},\tilde{\mathcal{G}})  $ defines a decomposition for the orbit segments on $M^f\times\mathbb{N}$.

	\subsection{\texorpdfstring{Specification on the good collection $ \tilde{\mathcal{G}} $}{Specification on the good collection $ \tilde{\mathcal{G}} $}}
	\begin{prop}\label{4.3}
		$ \tilde{\mathcal{G}} $ has specification at scale $ \hat{\delta} $ for any  $ \hat{\delta} $ sufficiently small.
	\end{prop}
	In order to prove Proposition 4.3, we will show that $ \mathcal{G} $ also has specification first. By Theorem IV.2.3 of \cite{Qian} and \cite{Cantarino} we know that there exists the local product structure, i.e. let $ \kappa>0 $ be sufficiently small, there is a $ \delta=\delta(\kappa)\in(0,\kappa) $ such that for each $ \tilde{x},\tilde{y}\in M^f $, if $ d(x,y)<\delta $ with $ \Pi(\tilde{x})=x, \Pi(\tilde{y})=y$, then
	$$W^u_{\kappa}(\tilde{x})\cap V^{cs}_\kappa(\tilde{y}) $$
	is a singleton, where $ V^{cs}_\kappa(\tilde{y})=\cup_{z\in \gamma^c_\kappa(\tilde{y})}W^s_\kappa(z) $ and  $ \gamma^c_{\kappa}(\tilde{y}) $ is a curve centered at $ y $ with radius less than $ \kappa $ and $ T_y\gamma^c_{\kappa}(\tilde{y})=E^c(\tilde{y}) $.
	Indeed, the existence of $\gamma^c_{\kappa}$ depends on the continuity of the center bundle. It can be seen that the disk $ V^{cs}_\kappa(\tilde{y}) $ is tangent to $ E^s\oplus E^c $ from Proposition 2.6 in \cite{BW}. In view of this,  let $\beta$ be small enough such that $\lambda\cdot \beta<\delta(\kappa)<\kappa$, where $\lambda=\max_{\tilde{x}\in M^f}\lVert Df|_{E^u(\tilde{x})} \rVert$ and
	we can choose some  curves $ \gamma^c_\beta(\tau^n\tilde{x}) $ which tangent to $ E^c(\tau^n\tilde{x}) $ and $f\gamma^c_\beta(\tau^{n-1}\tilde{x})\cap \gamma^c_\beta(\tau^n\tilde{x})$ contains an open interval of $\tau^n\tilde{x}$
	for every $ \tilde{x}\in M^f $.
	
	\begin{lem}
		Suppose that $ \tilde{\mathcal{W}}^u $ is $ \varepsilon$-minimal for any $ \varepsilon>0 $. Then for any $ \beta>0 $ small enough there exists $ L>0 $ such that for every $ \tilde{x}, \tilde{y}\in M^f $ if a disk $ D\subset W^u(\tilde{x}) $ with internal radius larger than $ L $, then $ D\cap V^{cs}_\beta(\tilde{y})\neq\emptyset $.
	\end{lem}
	
	\begin{pf}
		Fix $ \beta>0 $, let $ \delta=\delta(\beta) $ be sufficiently small such that for each $ \tilde{x},\tilde{y}\in M^f $, if $ d(x,y)<\delta $ we have
		$ W^u_\beta(\tilde{x})\cap V^{cs}_\beta(\tilde{y})\neq\emptyset $. We can get $ \mathcal{W}^u $ is $ \varepsilon$-minimal from the unstable foliation $ \tilde{\mathcal{W}}^u $ is $ \varepsilon$-minimal.
		For any $ \varepsilon<\delta $, there exist $L'>0$ such that if a disk $ D'\subset W^u(\tilde{x}) $ with internal radius larger than $ L' $ then
		$ D'\cap B(y,\varepsilon)\neq\emptyset $ for every $\tilde{x}\in M^f, y\in M $. Then for each $ \tilde{x},\tilde{y}\in M^f $ with  $ \Pi(\tilde{\sigma})=\sigma , \sigma\in\{x,y\}$,  there is a point $ \tilde{z}\in M^f $ so that $ z=\Pi(\tilde{z})\in D'\cap B(y,\varepsilon)  $ and $ W^u_\beta(\tilde{z})\cap  V^{cs}_\beta(\tilde{y})\neq\emptyset $. Thus,
		there exists a disk $ D\subset W^u(\tilde{x}) $ with internal radius larger than $ L $ and	$  D\cap V^{cs}_\beta(\tilde{y})\neq\emptyset $, where $ L=L'+\beta $.
	\end{pf}
	
	Before using the contraction property on the center direction for the good orbit segments in $(M^f,\tau)$, we define a corresponding good collection in $M$ as
	$$ \mathcal{G}:=\{(x,n)\in M\times\mathbb{N}\;|\; (\tilde{x},n)\in\tilde{\mathcal{G}} \text{ with } \Pi(\tilde{x})=x\},$$
	and then proceed to prove that the specification property holds on $\mathcal{G}$.
	\begin{lem}
		Fix $ r>0 $ as in Lemma 4.2, there exists $ \delta_0=\delta_0(r)>0 $ such that if $ (x,n)\in \mathcal{G} $, then $ V^{cs}_\delta(\tilde{x})\subset B_n(x,2\delta)$ for every $(\tilde{x},n)\in\tilde{\mathcal{G}}$ with $\Pi(\tilde{x})=x$ and all $ \delta<\frac{\delta_0}{2} $.
	\end{lem}
	\begin{pf}
		Fix $ x\in M $, as we know that Lemma 2.5 in \cite{Costa} implies the center-stable direction $ E^{cs}=E^s\oplus E^c $ only depends on $ x $ not $ \tilde{x} $ with $ \Pi(\tilde{x})=x $. For $ r>0 $, by the uniform continuity of $ \log\lVert Df|_{E^{cs}(x)}\rVert $ since $ M $ is compact, we can choose $ \delta_0>0 $ such that
		$$ \lvert \log\lVert Df|_{E^{cs}(x)}\rVert-\log\lVert Df|_{E^{cs}(y)}\rVert \rvert<\frac{r}{2} ,$$
		whenever $ d(x,y)<\delta_0 $ for every $ x,y\in M $.
		
		Let $ (x,n)\in \mathcal{G} $ and $ d_{V^{cs}(\tilde{x})} $ be the intrinsic distance defined by the Riemannian metric on $ V^{cs}_{\delta_0}(\tilde{x}) $. Since dim $E^c=1 $ and there is $ \tilde{x}\in M^f $ with $ \Pi(\tilde{x})=x $ such that $ S_j\tilde{\varphi}^c(\tilde{x})<-rj $ for all $ 1\leq j \leq n $, we have $ \lVert Df^j|_{E^{cs}(x)}\rVert \leq\lVert Df^j|_{E^{c}(\tilde{x})}\rVert<e^{-rj} ,\forall 1\leq j\leq n$ so that for $ \delta<\frac{\delta_0}{2} $ and $ y\in V^{cs}_\delta(\tilde{x}) $, then $ d(x,y)<2\delta $ and $ \lVert Df|_{E^{cs}(y)} \rVert<e^{-\frac{r}{2}} $. Thus, we have
		\begin{align*}
			d(fx,fy)\leq d_{V^{cs}(\tau\tilde{x})}(fx,fy)\leq e^{-\frac{r}{2}}d_{V^{cs}(\tilde{x})}(x,y)<2\delta
		\end{align*}
		and so that $ \lVert Df^2|_{E^{cs}(y)}\rVert\leq  e^{r}\lVert Df|_{E^{cs}(x)}\rVert\lVert Df|_{E^{cs}(fx)}\rVert<e^{-r} $. And so on we get $ \lVert Df^j|_{E^{cs}(y)}\rVert\leq e^{-\frac{r}{2}j} $ for all $ 1\leq j\leq n $, then $ d(f^{j}x,f^{j}y)\leq d_{V^{cs}(\tau^{j}\tilde{x})}(f^{j}x,f^{j}y)\leq e^{-\frac{r}{2}j}d_{V^{cs}(\tilde{x})}(x,y)<2\delta $, i.e. $ V^{cs}_\delta(\tilde{x})\subset B_n(x,2\delta)$.
		This completes the proof of this lemma.
	\end{pf}
	
	\begin{lem}
		$ \mathcal{G} $ has specification at scale $ 4\delta $ for any  $ 0<\delta<\delta_0/2 $ sufficiently small where $ \delta_0 $ is the constant defined in Lemma 4.5.
	\end{lem}
	
	\begin{pf}
		For any $ 0<\delta<\delta_0/2$, choose $ \frac{L}{2}>0 $ which was found in Lemma 4.4. Since $ f $ is uniformly expanding along the unstable manifolds  there is $ T\in \mathbb{N} $ such that for every $ \tilde{x}\in M^f $ and  $ n\geq T $, we have
		$ f^n W^u_\delta(\tilde{x})\supset D $ where $ D \subset W^u(\tau^n\tilde{x})$ is a disk with internal  radius larger than $ L/2 $ but less than $ L $. Then Lemma 4.4 implies that we can get  $ f^n W^u_\delta(\tilde{x})\cap V^{cs}_\delta(\tilde{y})\neq \emptyset $ for every $ \tilde{x},\tilde{y}\in M^f $ and $n\geq T$. Therefore, we choose a constant $ T=T(\delta)\in \mathbb{N} $ such that for every $ \tilde{x},\tilde{y}\in M^f $,
		\begin{equation}
			f^T W^u_\delta(\tilde{x})\cap V^{cs}_\delta(\tilde{y})\neq \emptyset
		\end{equation}
		and
		$$
		d(y_{-T},z_{-T})<\frac{1}{2}d(y_0,z_0)  \text{ whenever }  y_0, z_0\in W^u_\delta(\tilde{x}),
		$$
		where $ \tilde{\sigma}=(\sigma_n)_{n\in \mathbb{Z}} , \sigma\in\{y,z\}$ satisfying $ \Pi(\tilde{\sigma})=\sigma_0 , d(y_{-n},z_{-n})<2\delta$ and $ \lim\limits_{n\to \infty}d(y_{-n},z_{-n})=0 $.
		
		Now given $ \{(x^i,n_i)\}_{i=1}^k\subset \mathcal{G}  $,  then there exists a sequence $\{\tilde{x^i},n_i)\}_{i=1}^k \subset \tilde{\mathcal{G}}$.
		By (4) and the definition of the unstable manifold about a partially hyperbolic endomorphism, we can choose $ z^i_0, 1\leq i\leq k $ as follows. Set
		\begin{align*}
			\tilde{z}^1=\tilde{x}^1,\\
			z^2_0\in W^u_\delta(\tau^{n_1}\tilde{z}^1), z^2_{T}=f^{T}z^2_0\in V^{cs}_\delta(\tilde{x}^2) &\text{ with }  \tilde{z}^2 \text{ satisfying } d(z^2_{-n},z^1_{n_1-n})<\delta, n\geq 0,\\
			z^3_0\in W^u_\delta(\tau^{n_2+T}\tilde{z}^2)
			, z^3_{T}=f^{T}z^3_0\in V^{cs}_\delta(\tilde{x}^3) &\text{ with }  \tilde{z}^3 \text{ satisfying } d(z^3_{-n},z^2_{n_2+T-n})<\delta, n\geq 0,\\
			&\vdots \\
			z^k_0\in W^u_\delta(\tau^{n_{k-1}+T}\tilde{z}^{k-1}),
			z^k_{T}=f^{T}z^k_0\in V^{cs}_\delta(\tilde{x}^k) &\text{ with }\tilde{z}^k \text{ satisfying } d(z^k_{-n},z^{k-1}_{n_{k-1}+T-n})<\delta, n\geq 0.
		\end{align*}
		Then  $ \tilde{z}^k $ could shadow the orbit segments $ \{(x^i,n_i)\}_{i=1}^k $ we gave before in the sense of the corresponding time. Indeed, according to the choice of $ T $, we have if $ z^i_0\in W^u_\delta(\tau^{n_{i-1}+T}\tilde{z}^{i-1}) ,(3\leq i\leq k)$, then
		\begin{gather*}
			d(z^i_{-j},z^{i-1}_{n_{i-1}+T-j})<\delta, 1\leq j\leq n_{i-1},\\
			d(z^i_{-n_{i-1}-T-j},z^{i-1}_{-j})<\frac{\delta}{2}, 1\leq j\leq n_{i-2},\\
			\vdots\\
			d(z^i_{-n_2-\cdots-n_{i-1}-(i-2)T-j},z^{i-1}_{-n_2-n_3-\cdots-n_{i-2}-(i-3)T-j})<\frac{\delta}{2^{i-2}}, 1\leq j\leq n_{1}.
		\end{gather*}
		Thus, together with $ \{(x^i,n_i)\}_{i=1}^k\subset \mathcal{G} $ has uniform contraction in the direction $ E^{cs} $, we have
		$$
		d(f^{j}z^k_{-m_{k-i}},f^j x^i)<2\delta+\sum_{i=0}^{\infty}\frac{\delta}{2^i}=4\delta, 0\leq j\leq n_i-1 , 1\leq i\leq k-1,
		$$
		and
		$ d(f^{j+T}(z^k_0),f^j x^k)<4\delta, 0\leq j\leq n_k-1 $,
		where $ -m_{k-i}=-\sum_{j=i}^{k-1}n_j-(k-i-1)T $,
		i.e. $ \mathcal{G}  $ has specification at scale $ 4\delta $. This completes the proof of this lemma.
	\end{pf}

	\begin{pot4.3}
		For any $ \hat{\delta}>0 $ there exist an integer $ N\geq 0 $ and $ \delta_1>0 $ such that if $ \tilde{x}=(x_n)_{n\in \mathbb{N}} , \tilde{y}=(y_n)_{n\in \mathbb{N}}\in M^f $ and $ d(x_0,y_0)<\delta_1 $, then $ \tilde{d}(\tau^N \tilde{x}, \tau^N \tilde{y})<\hat{\delta} $.
		Given $ \{(\tilde{x}^i,n_i)\}_{i=1}^k\subset \tilde{\mathcal{G}} $ where $ \tilde{x}^i=(x^i_n)_{n\in\mathbb{Z}} $ with $ x^i_0=x^i $ then $ \{(x^i,n_i)\}\subset \mathcal{G} $, next we will show that  $ \{(x^i_{-N},n_i+N)\}_{i=1}^k $ has specification similarly as Lemma 4.6.
		
		Without loss of generality, we assume that the larger rate $ \lambda_2 $ on $ E^{cs} $ is greater than or equal to 1. If not, the orbits on $ E^{cs} $ are uniformly contracting, which is the same as Lemma 4.6. Let $ \delta<\min\{\delta_0/2,\delta_1/4\} $ where $ \delta_0 $ defined in Lemma 4.5. For any $ \tilde
		x \in M^f$, we take $ \gamma^{c}_{\lambda_2^{-j}\delta}(\tau^{-j}\tilde{x}) $ satisfying $ f\gamma^{c}_{\lambda_2^{-j}\delta}(\tau^{-j}\tilde{x})\cap \gamma^{c}_{\lambda_2^{-j+1}\delta}(\tau^{-j+1}\tilde{x}) $ contains an open interval of $ x_{-j+1}=\Pi(\tau^{-j+1}\tilde{x}) (1\leq j\leq N) $ since $ D_{x_n}fE^c(\tau^n \tilde{x})=E^c(\tau^{n+1}\tilde{x}) $ for any $ n\in\mathbb{Z} $. Let $ \gamma^c_\delta(\tilde{x})\subset V^{cs}_\delta(\tilde{x}) $ and  define $ V_{-N}(\tilde{x}):=\bigcup_{z\in \gamma^c_{\lambda_2^{-N}\delta}(\tau^{-N}\tilde{x})}W^s_{\delta}(z) $.
		For $\lambda^{-N}_2\delta>0$ there exists $\delta'\in (0,\lambda^{-N}_2\delta)$ such that there is local product structure whenever $d(x,y)<\delta'$.
		Suppose that  $ \mathcal{W}^u $ is $ \alpha$-minimal where $ \alpha<\delta' $. Then we can choose a  integer $ T=T(\delta,N)\in \mathbb{N} $ large enough such that
		\begin{align*}
			f^TW^u_\delta(\tilde{x})\cap V_{-N}(\tilde{y})\neq \emptyset  \text{ for every }   \tilde{x},\tilde{y}\in M^f
		\end{align*}
		and
		$$
		d(y_{-T},z_{-T})<\frac{1}{2}d(y_0,z_0)  \text{ whenever }  y_0, z_0\in W^u_\delta(\tilde{x}),
		$$
		where $ \tilde{\sigma}=(\sigma_n)_{n\in \mathbb{Z}} , \sigma\in\{y,z\}$ satisfying $ \Pi(\tilde{\sigma})=\sigma_0 , d(y_{-n},z_{-n})<2\delta$ and $ \lim\limits_{n\to \infty}d(y_{-n},z_{-n})=0 $.
		Next we construct $ \tilde{y}^1, \tilde{y}^2,\cdots, \tilde{y}^k $ so that the orbit segments of $ \tilde{y}^k $ could shadow $ \{(x^i_{-N},n_i+N)\}_{i=1}^k$. We pick up  $\tilde{y}^1, \tilde{y}^2,\cdots, \tilde{y}^k  $ in such way
		\begin{align*}
			\tilde{y}^1=\tilde{x}^1,\\
			y^2_0\in W^u_\delta(\tau^{n_1}\tilde{y}^1), y^2_{T}=f^{T}y^2_0\in V_{-N}(\tilde{x}^2) &\text{ with } \tilde{y}^2 \text{ satisfying } d(y^2_{-n},y^1_{n_1-n})<\delta, n\geq 0,\\
			y^3_0\in W^u_\delta(\tau^{n_2+T+N}\tilde{y}^2), y^3_{T}=f^{T}y^3_0\in V_{-N}(\tilde{x}^3) &\text{ with }  \tilde{y}^3 \text{ satisfying } d(y^3_{-n},y^2_{n_2+T+N-n})<\delta, n\geq 0,\\
			&\vdots \\
			y^k_0\in W^u_\delta(\tau^{n_{k-1}+T+N}\tilde{y}^{k-1}), y^k_{T}=f^{T}y^k_0\in V_{-N}(\tilde{x}^k) &\text{ with } \tilde{y}^k \text{ satisfying } d(y^k_{-n},y^{k-1}_{n_{k-1}+T+N-n})<\delta, n\geq 0.
		\end{align*}
		For every $ \omega\in V_{-N}(\tilde{x}^i) $,
		$ d(f^j\omega, f^j x_{-N+j}^i)<2\delta $ for all $ 0\leq j\leq N $ and $ f^N\omega \in V^{cs}_\delta(\tilde{x}^i) $. Moreover, $ \{(x^i,n_i)\}\subset \mathcal{G} $, so $ f^N(\omega)\in B_{n_i}(x^i,2\delta) $, i.e. $ \omega\in B_{n_i+N}(x^i_{-N},2\delta) $. Set $ -m_{k-i}=-\sum_{j=i}^{k-1}n_{j}-(k-i)N-(k-i-1)T$ for any $ 1\leq i \leq k-1$,
		then we have
		$$
		d(f^{j}y^k_{-m_{k-i}},f^j x^i_{-N})<2\delta+\sum_{i=0}^{\infty}\frac{\delta}{2^i}=4\delta, 0\leq j\leq n_i-1+N , 1\leq i\leq k-1,
		$$
		and
		$$
		d(f^{j+T}y^k_0,f^j x^k_{-N})<4\delta, 0\leq j\leq n_k-1.
		$$
		Therefore, use the fact we mentioned at first, we get
		$$
		\tilde{d}(\tau^{-m_{k-i}+N+j}\tilde{y}^k,\tau^j \tilde{x}^i)<\hat{\delta}, 0\leq j \leq n_i-1, 1\leq i\leq k-1,
		$$
		and
		$$
		\tilde{d}(\tau^{j+T+N}\tilde{y}^k,\tau^j \tilde{x}^k)<\hat{\delta}, 0\leq j\leq n_k-1.
		$$
		This completes the proof of this proposition.
	\end{pot4.3}

	\subsection{\texorpdfstring{The pressure of obstructions to expansivity}{The pressure of obstructions to expansivity}}
	It is well known that the bi-infinite Bowen ball around $x$  has to be contained in a compact subset of the center leaf of $ x $ when $ f $ is a partially hyperbolic diffeomorphism with dim$ (E^c)=1 $ \cite{LVY}. For a partially hyperbolic endomorphism, we need to use the condition that $ f $ is non-degenerate to control the preimages.
	
	Let $ \varepsilon_0 $ be an exponent of separation for $ f $. For any $ \varepsilon>0 $ and $\tilde{x}$ with $ \Pi(\tilde{x})=x_0 $, define a set
	$$
	\Gamma_{\varepsilon}(\tilde{x})=\{y_0\in M\;|\;\tilde{y}\in M^f \text{ with } \Pi(\tilde{y})=y_0 \text{ such that } d(x_n,y_n)<\varepsilon, \forall n\in \mathbb{Z}\}.
	$$
	
	\begin{lem}
		Given an ergodic measure $ \mu\in \mathcal{M}^e_f $, for every $ 0<r<-\log\lambda_s $ there is $ \varepsilon(r)>0 $ such that if for $ \mu $-$a.e.$ $ x_0\in M $,
		$ \lambda^c(\mu,f)=\lim\limits_{n\to \infty}\frac{1}{n}\log \lVert D_{x_0}f^n|_{E^c(\tilde{x})}\rVert<-r$, then
		$\Gamma_{\varepsilon(r)}(\tilde{x})=\{x_0\}$, $ \mu $-$a.e.$ $ x_0\in M $ with $ \Pi(\tilde{x})=x_0 $.
	\end{lem}
	\begin{pf}
		Firstly, we show that for any $ \beta>0$ sufficiently small, there exists $ \varepsilon=\varepsilon(\beta)>0 $ such that $ \Gamma_{\varepsilon}(\tilde{x})\subset \gamma^c_\beta(\tilde{x})$ for every $ \tilde{x}\in M^f $.
		Take $ 0<\varepsilon<\delta(\beta) $ appropriately, so that local product structure exists in $B(x,\delta(\beta))$ for any $x\in M$.
		For every $ \tilde{x}\in M^f $ and $ y_0\in \Gamma_{\varepsilon}(\tilde{x}) $ there exists $ \tilde{y}\in M^f $ with $ \Pi(\tilde{\sigma})=\sigma_0, \sigma\in\{x,y\} $ such that
		$$
		d(x_n,y_n)<\varepsilon\text{ for all }n\in\mathbb{Z},
		$$
		and there are $ \tilde{z}^i=(z^i_n)_{n\in\mathbb{Z}}  \in M^f$ $(i=1,2)$ with
		$$
		z^1_n\in V^{cs}_\beta(\tau^n\tilde{x})\cap W^u_\beta(\tau^n\tilde{y}), \;
		z^2_n\in W^s_\beta(z^1_n)\cap\gamma^c_\beta(\tau^n\tilde{x}).
		$$
		Under forward iterates, the distance between $ y_j $ and $ z^1_j $ is expanding and
		\begin{align*}
			d(y_j, z^1_j)
			&<\varepsilon+d(x_j, z^2_j)+d(z^2_j, z^1_j),
		\end{align*}
		if $ y_0\neq z^1_0 $ there must exist some $ n\geq 0 $ is large enough that contradicts the inequality above, since the distance between $ x_j $ and $ z^2_j $ is controlled and the distance between $z^2_j$ and $z^1_j$ is decreasing for any $j\geq 0$. Similarly, we consider under backward iterates because there is uniform expansion along $ W^s $. We suppose that $ y_0\neq z^2_0 $, then for any $n\geq 0$,
		\begin{equation}
			d(y_{-n},z^2_{-n})\leq
			d(y_{-n},x_{-n})+d(x_{-n},z^2_{-n})<\varepsilon+\beta
		\end{equation}
		However, the point  $ y_{-n}\in W^s(z^2_{-n}) $ for any $ n\geq 0 $.
		Thus, the distance between $ y_{-n} $ and $ z^2_{-n} $ could expand which contradicts with (5).
		Therefore, we get $ \Gamma_{\varepsilon}(\tilde{x})\subset \gamma^c_{\beta}(\tilde{x}) $ for every $ \tilde{x}\in M^f $.
		
		Now we consider the condition about the central Lyapunov exponent. For a fixed  $0<r<-\log\lambda_s  $,
		from Proposition I.3.5 in \cite{Qian} we know that for $ \mu$-$a.e.$ $ x_0\in M $,
		$ \lambda^c(\tilde{x})=\lambda^c(x_0)<-r $ for every $\tilde{x}\in \Pi^{-1}(x_0)$, then  we consider the exponential rate of  backward iterates about $ D_{x_0}f|_{E^c(\tilde{x})} $ whenever given $ \tilde{x} $. Since the central Lyapunov exponent at $ \tilde{x} $ is given by
		\begin{equation*}
			\lambda^c(\tilde{x})=\lim\limits_{m\to-\infty}\frac{1}{m}\log\lVert D_n^m(\tilde{x})|_{E^c(\tau^n\tilde{x})}\rVert,
		\end{equation*}
		where $ D_n^m(\tilde{x})=(D^{-m}_{n+m})^{-1}=(D_{x_{n+m}}f)^{-1}|_{E^c(\tau^{n+m+1}\tilde{x})}\circ \cdots\circ (D_{x_{n-1}}f)^{-1}|_{E^c(\tau^n\tilde{x})}$ if $ m<0 $, then we have
		\begin{equation*}
			-\lambda^c(\tilde{x})=\lim\limits_{m\to+\infty}\frac{1}{m}\log\lVert D_{x_{n}}f^{-m}|_{E^c(\tau^{n}\tilde{x})}\rVert>r
		\end{equation*}
		for any $ n\in\mathbb{Z} $. For $ n=0$, there are arbitrarily large $ m $ and some constant $ C>0 $ such that
		\begin{equation}
			\lVert D_{x_0}f^{-m}|_{E^c{(\tilde{x})}} \rVert>Ce^{rm}.
		\end{equation}
		For any $\beta>0$ sufficiently small, choosing $ 0<\varepsilon(r)<\min\{\varepsilon_0/6,\varepsilon=\varepsilon(\beta)\} $  so that $ \Gamma_{\varepsilon(r)}(\tilde{x})\subset \gamma^c_\beta(\tilde{x}) $ and $ f|_{B(x,3\varepsilon(r))}:B(x,3\varepsilon(r))\to f(B(x,3\varepsilon(r))) $ is homeomorphism as $ \varepsilon_0 $ is the exponent of separation for $ f $. Without loss of generality,  we suppose that $ \varepsilon(r) $ is sufficiently small such that
		$$ \lvert \log\lVert D_{x_0}f^{-1}|_{E^c{(\tilde{x})}}\rVert-  \log\lVert D_{x'_0}f^{-1}|_{E^c{(\tilde{x}')}}\rVert \rvert <\frac{r}{2}$$
		whenever $ \tilde{d}(\tau^j\tilde{x},\tau^j\tilde{x}')<3\varepsilon(r) $ for all $ j\in \mathbb{Z} $ since the space $ E^c(\tilde{x}) $ varies continuously. Then for any $ y_0\in \Gamma_{\varepsilon(r)}(\tilde{x}) $, that is, $ \tilde{d}(\tau ^j\tilde{x},\tau^j\tilde{ y})=\tilde{d}(\tilde{x},\tilde{y})<3\varepsilon(r) $, so  for the $ m $ chosen satisfying (6) we have
		$$
		\lVert D_{y_0}f^{-m}|_{E^c{(\tilde{y})}}\rVert>\lVert D_{x_0}f^{-m}|_{E^c{(\tilde{x})}}\rVert\cdot e^{-\frac{r}{2}m}>Ce^{\frac{r}{2}m}.
		$$
		Since $ \Gamma_{\varepsilon(r)}(\tilde{x})\subset \gamma^c_\beta(\tilde{x}) $, then
		$$
		d(x_{-m},y_{-m})\geq Ce^{\frac{r}{2}m}d(x_0,y_0)
		$$
		for all $ m>0 $ arbitrarily large satisfying (6). Therefore, $ x_0=y_0 $, i.e. $ \Gamma_{\varepsilon(r)}(\tilde{x})=\{x_0\} $ for $ \mu $-$a.e.$ $ x_0\in M $. This completes the proof of this lemma.
	\end{pf}
	Let $\tilde{\mu}$ be an ergodic measure with $ \Pi(\tilde{\mu})=\mu $. From the proof above, we know $ \varepsilon(r)<\varepsilon_0 $ for $ \varepsilon(r) $ in Lemma 4.7 and $ \varepsilon_0 $ is the separation exponent for $ f $. We claim that $ \tilde{\mu} $ is almost $\tau$-expansive at scale $ \varepsilon(r) $ under the condition of Lemma 4.7. This is because that for $ \tilde{\mu} $-a.e. $ \tilde{x}=(x_n)_{n\in\mathbb{Z}}, \tilde{y}=(y_n)_{n\in\mathbb{Z}}\in M^f $ if $\tilde{d}(\tau^n \tilde{x},\tau^n \tilde{y})<\varepsilon(r) $ for all $ n\in \mathbb{Z} $ then by the definition of metric $ \tilde{d} $ on $ M^f $ we have $ d(x_n,y_n)<\varepsilon(r) $ for all $ n\in \mathbb{Z} $, so $ y_0\in \Gamma_{\varepsilon(r)}(\tilde{x}) $. It follows that there must have $ x_0=y_0 $ and as $ \varepsilon(r)<\varepsilon_0 $ we get $ x_{-n}=y_{-n} $ when $ n\geq 0 $ so that $ \tilde{x}=\tilde{y} $.
	
	For any $\varepsilon>0$ and every $ \tilde{x}\in M^f $, let us define the set $ \tilde{\Gamma}_{\varepsilon}(\tilde{x})=\{\tilde{y}\in M^f \;|\;\tilde{d}(\tau^n \tilde{x},\tau^n \tilde{y})<\varepsilon, n\in \mathbb{Z}\} $ and the non-expansive set of $ \tau $ at scale $ \varepsilon $ by
	$ \tilde{NE}(\varepsilon)=\{\tilde{x}\in M^f \;|\;\tilde{\Gamma}_{\varepsilon}(\tilde{x})\neq \{\tilde{x}\} \} $. We know that Lemma 4.7 implies the set of ergodic measures on $ M^f $ which are not almost $ \tau$-expansive at scale $ \varepsilon(r) $ is contained in the set of ergodic measures with central Lyapunov exponent not less than $ -r $. Therefore, by Lemma 4.2 there exists $ \varepsilon(r)>0 $ such that
	\begin{align*}
		P^\perp_{\mathrm{exp}}(\tau, \phi^*, \varepsilon(r))&=\sup_{\tilde{\mu}\in \mathcal{M}^e_{\tau}}\{h_{\tilde{\mu}}(\tau)+\int \phi^* d\tilde{\mu} \;|\; \tilde{\mu}(\tilde{NE}(\varepsilon(r)))>0\}\\
		&\leq\sup_{\tilde{\mu}\in \mathcal{M}^e_{\tau}}\{h_{\tilde{\mu}}(\tau)+\int \phi^*d\tilde{\mu} \;|\; \lambda^c(\tilde{\mu},\tau)\geq -r\}\\
		&<P(\tau,\phi^*).
	\end{align*}

	\subsection{\texorpdfstring{The pressure of obstructions to specification}{The pressure of obstructions to specification}}
	In order to prove that there exists a unique equilibrium state for $ (M^f,\tau,\phi^*)$, we will show that the obstructions to specification $ \tilde{\mathcal{P}} $ have a small pressure in this subsection.
	\begin{prop}\label{prei}
		There exists $ \varepsilon_1>0 $ such that $ P(\tilde{\mathcal{P}}, \phi^*,\delta,\varepsilon)<P(\tau,\phi^*) $ for every $ 0<\varepsilon<\varepsilon_1 $ and any $\delta>0$.
	\end{prop}
	
	\begin{pf}
		For $\delta>0$, let $ E_n\subset \tilde{\mathcal{P}}_n$ be an $(n,\delta) $-separated set with
		\begin{equation*}
			\log \sum_{\tilde{x}\in E_n}e^{\Phi^*_0(\tilde{x},n)}\geq \log \Lambda(\tilde{\mathcal{P}},\tau,\phi^*,\delta,n)-1
		\end{equation*}
		and $ \tilde{\sigma}_n\in \mathcal{M}(M^f) $ be the atomic measure concentrated on $ E_n $ by the formula
		\begin{equation*}
			\tilde{\sigma}_n=\frac{ \sum_{\tilde{x}\in E_n}e^{\Phi^*_0(\tilde{x},n)}\delta_{\tilde{x}}}{\sum_{\tilde{z}\in E_n}e^{\Phi^*_0(\tilde{z},n)}},
		\end{equation*}
		where $\Phi^*_0(\tilde{x},n)=\sum_{i=0}^{n-1}\phi^*(\tau^i\tilde{x})$.
		Let $ \tilde{\mu}_n\in \mathcal{M}(M^f) $  be defined by $ \tilde{\mu}_n=\frac{1}{n}\sum_{i=0}^{n-1}\tilde{\sigma}_n\circ \tau^{-i} $. By  Theorem 9.10 of \cite{Walters}, we obtain that any limit point $ \tilde{\mu} $ of $\{\tilde{\mu}_n\} $ is $ \tau$-invariant and has
		$$
		h_{\tilde{\mu}}(\tau)+\int \phi^* d\tilde{\mu}
		\geq
		P(\tilde{\mathcal{P}},\phi^*,\delta).
		$$
		Set $ \mu:=\Pi(\tilde{\mu}), \mu_n:=\Pi(\tilde{\mu}_n) $. For every $ (\tilde{x},n)\in E_n\subset \tilde{\mathcal{P}} $,  that is, $ S_n\tilde{\varphi} ^c(\tilde{x})\geq -rn $, then we get
		\begin{align*}
			\int \tilde{\varphi}^c d\tilde{\mu}_n&=\frac{1}{n}\sum_{i=0}^{n-1}\int \tilde{\varphi}^c \circ \tau^i d\tilde{\sigma}_n\\
			&=\frac{1}{n}\sum_{i=0}^{n-1} \tilde{\varphi}^c\circ \tau^i(\tilde{x}), \; \tilde{x}\in E_n\\
			&=\frac{1}{n} S_n\tilde{\varphi}^c(\tilde{x})\geq -r,
		\end{align*}
		so $ \lambda^c(\tilde{\mu}_n,\tau)\geq -r $ and by the weak*-convergence we have $ \lambda^c(\tilde{\mu},\tau)\geq -r $. Then by Lemma 4.2, we get
		\begin{align*}
			P(\tau,\phi^*)>h_{\tilde{\mu}}(\tau)+\int \phi^* d\tilde{\mu}\geq
			P(\tilde{\mathcal{P}},\phi^*,\delta).
		\end{align*}
		Therefore, there exists $ \varepsilon_1>0 $ such that for any $ 0<\varepsilon<\varepsilon_1 $, we have
		$ P(\tau,\phi^*)>P(\tilde{\mathcal{P}},\phi^*,\delta,\varepsilon) $. This completes the proof of this proposition.
	\end{pf}

	\subsection{\texorpdfstring{The Bowen property on $\tilde{\mathcal{G}}$}{The Bowen property on $\tilde{\mathcal{G}}$}}
	\begin{prop}
		For every H\"{o}lder continuous function $ \phi: M\to \mathbb{R} $, there exists $ \tilde{\varepsilon}>0 $ such that $ \phi^* $ has Bowen property on $ \tilde{\mathcal{G}} $ at any scale $ 0<\varepsilon<\tilde{\varepsilon}$.
	\end{prop}
	
	\begin{pf}
		For any $\beta'<\min\{\beta, \delta_0/2\}$ where $\beta$ is defined in §4.3 and $\delta_0$ as in Lemma 4.5, there exists $\tilde{\varepsilon}:=\delta(\beta')$ such that if $d(x_0,y_0)<\tilde{\varepsilon}$, then
		\begin{align*}
			V^{cs}_{\beta'}(\tilde{x})\cap W^u_{\beta'}(\tilde{y})
		\end{align*}
		is a singleton for every $\tilde{x},\tilde{y}\in M^f$ with $\tilde{\sigma}=(\sigma_n)$, $\sigma\in\{x,y\}$.
		For any $\varepsilon<\tilde{\varepsilon}$, given $(\tilde{x},n)\in \tilde{\mathcal{G}}$, for any $\tilde{y}\in B_n(\tilde{x},\varepsilon)$, that is, $\tilde{d}(\tau^i\tilde{x},\tau^i\tilde{y})<\varepsilon$ for any $0\leq i\leq n-1$. There exist $z_i\in V^{cs}_{\beta'}(\tau^i\tilde{x})\cap W^u_{\beta'}(\tau^i\tilde{y})$ satisfying $f(z_i)=z_{i+1}$, $0\leq i\leq n-2$. By Lemma 4.5, we have
		\begin{align*}
			d(x_i,z_i)\leq e^{-\frac{r}{2}i}\cdot 2\beta' \text{ for all }0\leq i\leq n-1.
		\end{align*}
		Since $ z_{n-1}\in W^u_{\beta'}(\tau^{n-1}\tilde{y}) $ then we have
		$$
		d(z_i,y_i)\leq \lambda^{-(n-i-1)}_ud(z_{n-1},y_{n-1})\leq \lambda^{-(n-i-1)}_u\beta',\; i\in\{0,1,\cdots,n-1\}.
		$$
		Then the triangle inequality gives
		$$
		d(x_i,y_i)\leq 2\max\{e^{-\frac{r}{2}i}\cdot 2\beta', \lambda^{-(n-i-1)}_u\beta'\}
		$$
		for every $ 0\leq i\leq n-1 $. We note that $ S_n\phi^*(\tilde{x})=S_n\phi(x_0) $ and let $ K $ be the H\"{o}lder constant of $ \phi $ and $ \alpha\in (0,1) $ be the H\"{o}lder exponent, we can have
		\begin{align*}
			\lvert S_n\phi^*(\tilde{x})-S_n\phi^*(\tilde{y})\rvert
			&=\lvert S_n\phi(x_0)-S_n\phi(y_0)\rvert\\
			&\leq K\sum_{i=0}^{n-1} (4\beta')^\alpha\max\{e^{-\frac{r}{2}i\alpha},\lambda^{-(n-i-1)\alpha}_u\}\\
			&\leq  K (4\beta')^\alpha\sum_{i=0}^{\infty}(e^{-\frac{r}{2}i\alpha}+\lambda^{-i\alpha}_u):=K'<\infty
		\end{align*}
		We can see that $ K' $ is independent of $ \tilde{x},\tilde{y} $ and $ n $. This completes the proof of this proposition.
	\end{pf}
	
	\subsection{\texorpdfstring{Proof of Main Theorem}{Proof of Main Theorem}}
	Take the map $f:M\to M$ and the potential function $\phi:M\to \mathbb{R}$ as in \thmref{Main}.
	Let $ \varepsilon<\min\{\varepsilon(r) ,\varepsilon_1,\tilde{\varepsilon}\}$, where $\varepsilon(r),\varepsilon_1, \tilde{\varepsilon}  $ are as in Lemma 4.7, Proposition 4.8 and 4.9 respectively, we have
	\begin{gather*}
		P^\perp_{\mathrm{exp}}(\tau,\phi^*,\varepsilon)<P(\tau,\phi^*),\\
		P(\tilde{\mathcal{P}},\phi^*,\delta,\varepsilon)<P(\tau,\phi^*),
	\end{gather*}
	and $ \phi^* $ has Bowen property on $ \tilde{\mathcal{G}} $ at scale $ \varepsilon $. By Proposition 4.3 we know that $ \tilde{\mathcal{G}} $ has specification at scale $ \delta $ with $ \varepsilon\geq 2000\delta $. Thus by Theorem \ref{2.2} we obtain that $ (\tau,\phi^*) $ has a unique equilibrium state $ \tilde{\mu} $. Moreover, let $ \Pi(\tilde{\mu}):=\mu $ we have
	\begin{equation*}
		P(\tau,\phi^*)=h_{\tilde{\mu}}(\tau)+\int \phi^* d\tilde{\mu}=h_{\mu}(f)+\int \phi d\mu
	\end{equation*}
	and
	$$
	P(\tau,\phi^*)=P(f,\phi)
	$$
	Therefore, $ \mu $ is also a unique equilibrium state for $ (f,\phi) $.
	
	Next, we show the robustness of the uniqueness of the equilibrium states. By Proposition \ref{3.1},
	for any $ 0<a<\frac{P^u(f,\phi)-P^s(f,\phi)}{2} $
	there exists a $ C^1 $ neighborhood  $ \mathcal{U} $ of $ f $ such that
	$ P^u(g,\psi)>P^s(g,\psi)$
	for any H\"{o}lder continuous function $\psi\in U_{a}(\phi)$ and any $ g\in \mathrm{NDPHE}^2(M)\cap \mathcal{U} $.
	We can choose a $ C^1 $-neighborhood $  \mathcal{V} $ of $ f $ such that the unstable foliation $ \tilde{\mathcal{W}}^u_g $ for $ g $ is $ \alpha$-minimal for every $ g\in\mathcal{V} $, where $ \alpha $ is chosen in the proof of Proposition 4.3. Set $ \mathcal{U}_f=\mathcal{U}\cap \mathcal{V} $, then for any partially hyperbolic endomorphism $ g\in \mathrm{NDPHE}^2(M)\cap \mathcal{U}_f $ and any H\"{o}lder continuous function $\psi\in U_{a}(\phi)$ there exists a unique equilibrium state for $ (g,\psi) $. This completes the proof of \thmref{Main}.

	\section{Example}
	
	In this section, we will give an example which is derived from an Anosov endomorphism satisfying the conditions of \thmref{Main}. These kinds of examples for diffeomorphisms were constructed by  Mañé in \cite{M78}.
	
	We fix a matrix $ A\in GL(d,\mathbb{Z}) $ with det $  A>1 $ and all eigenvalues simple,  irrational, and positive and only one eigenvalue outside the unit circle. In order to facilitate the discussion, we assume $ d=3 $, $ \lambda_u>1>\lambda_c>\lambda_s>0$ and the corresponding eigenspaces $ G^{u,c,s} \subset \mathbb{R}^3$. For example, let a matrix be defined by
	\begin{align*}
		A=\begin{pmatrix}
			100& 1 & 0\\
			-100 & 0& 1\\
			3 & 0 & 0\\
		\end{pmatrix},
	\end{align*}
	it is easy to check that such a matrix $ A $ meets the above conditions.
	Let $ f_{A} : \mathbb{T}^3\to \mathbb{T}^3 $ be the hyperbolic toral endomorphism induced by $ A $ and  $ \mathcal{F}^{u,c,s} $ be the foliations of $ \mathbb{T}^3 $ by the leaves parallel to the eigenspaces $ G^{u,c,s} $, respectively. Since all eigenvalues are irrational, every leaf of $ \mathcal{F}^u, \mathcal{F}^c $ and $\mathcal{F}^s $ is dense in $ \mathbb{T}^3 $.
	
	Let fix $ \rho>\rho'>0 $ such that $\tau$ on $(\mathbb{T}^3)^{f_A}$ is expansive at scale $\rho$ and $ q $ be a fixed point of $ f_{A} $, now we construct a small perturbation of $ f_{A} $ to get a partially hyperbolic endomorphism $ f $ on $ \mathbb{T}^3 $. Specifically, we set $ f=f_{A} $ outside of $  B(q,\rho) $ .  Inside $  B(q,\rho) $, let the fixed point $ q $ undergo a pitchfork bifurcation in the direction of the center foliation $ \mathcal{F}^c $. In the same way, $ f $ satisfies the following: for every $ \tilde{x}\in (\mathbb{T}^3)^f $,
	\begin{enumerate}
		\item the center foliation $ W^c(\tilde{x}):=\mathcal{F}^c(\tilde{x}) $ is $ f $-invariant, and $ T_{\tilde
			{x}}W^c(\tilde{x})=E^c(\tilde{x}) $;
		\item the cones around $ G^u_{\tilde
			x}, G^s_{\tilde
			x} $ are invariant and uniformly expanding for $ D_{\tilde{x}}f, D_{\tilde{x}}f^{-1} $ respectively, and $ D_{\tilde{x}}f $-invariant distributions $ E^u(\tilde{x}), E^s(\tilde{x}) $ integrating to $ f $-invariant  $ W^u(\tilde{x}), W^s(\tilde{x}) $;
		\item we have $ \lVert Df|_{E^{cs}} \rVert\leq \lambda_c <1 $ outside of $ B(q,\rho') $;
		\item for every $ \tilde{x}\in (\mathbb{T}^3)^f $, det $ D_{\tilde{x}}f \neq 0 $.
	\end{enumerate}
	Therefore, we get a partially hyperbolic endomorphism $ f $ on $ \mathbb{T}^3 $ with $ T_{\tilde{x}}\mathbb{T}^3=E^s(\tilde{x})\oplus E^c(\tilde{x})\oplus E^u(\tilde{x})  $ for every $ \tilde{x}\in  (\mathbb{T}^3)^f $. See more details about the construction of the perturbation in \cite{M78}. We note that
	\begin{align*}
		\lambda^c(f)=\sup\{\lVert Df|_{E^c(\tilde{x})}\rVert\;|\; \tilde{x}\in (\mathbb{T}^3)^f \}>1
	\end{align*} 		
	since the center direction experiences expansion inside $ B(q,\rho') $.
	We might consider $f_A$ as a partially hyperbolic endomorphism on $\mathbb{T}^3$, then we have
	\begin{align*}
		h^u(f_A)>h^s(f_A).
	\end{align*}
	By Proposition \ref{3.1}, we could control the perturbation to be sufficiently small so that $f$ also satisfies that $h^u(f)>h^s(f)$. Moreover, the unstable foliation of $ f $ is minimal, then we can use \thmref{Main} to prove that there exists a $ C^1 $-neighborhood $ \mathcal{U}_f $ of $ f $ such that each non-degenerate $ C^2 $ endomorphism $ g\in \mathcal{U}_f $ has a unique measure of maximal entropy.

	\section*{Acknowledgments}

	\footnotesize{

		\textit{ E-mail address:} \href{mailto:\Emaila}{zhangyf@stu.xmu.edu.cn}

		\textit{ E-mail address:} \href{mailto:\Emailb}{yjzhu@xmu.edu.cn}}

\end{document}